\documentclass[11pt,reqno]{amsart}

\usepackage{amssymb,amsmath,amsthm,graphics,epsfig}
\usepackage{amsfonts}
\usepackage{graphicx,cite}
\usepackage{amssymb}
\usepackage[pagewise]{lineno}
\usepackage{mathrsfs}
\usepackage[colorlinks=true, pdfstartview=FitV, linkcolor=blue, citecolor=blue, urlcolor=blue]{hyperref}
\usepackage{booktabs}
\usepackage{comment}
\usepackage{graphicx}
\usepackage{tikz}
\usepackage{multicol}
\usepackage{cases}
\usepackage{booktabs}
\usepackage{subcaption}
\usepackage{orcidlink}

\usepackage{arydshln}
\usepackage[colorinlistoftodos]{todonotes}
\usepackage[normalem]{ulem}
\usepackage{color}
\usepackage{amsaddr} %direcciones adelante
\usepackage[margin=2.7cm]{geometry}

\usepackage{float}
\theoremstyle{definition}

\definecolor{ColorEdward}{rgb}{0.4,0.6,0.7}

\definecolor{ColorFrancisco}{rgb}{0.956,0.878,0.5}

\definecolor{ColorJosep}{rgb}{0.2,0.5,0.1}

\definecolor{ColorEdward}{rgb}{0.4,0.6,0.7}
\definecolor{ColorFrancisco}{rgb}{0.956,0.878,0.5}
\definecolor{ColorJosep}{rgb}{0.2,0.5,0.1}

\title[KAM stability of 14 Herculis]{The Orbital Architecture of 14 Herculis System as a KAM-stable Three-body Configuration}

\author[F. Crespo, J. F. Palacián, E. A. Turner and P. Yanguas]{{ Francisco Crespo \orcidlink{0000-0002-5930-8523}$^{1}$, Jesús Palacián \orcidlink{0000-0002-0974-6656}$^2$, Edward A. Turner \orcidlink{0000-0002-2959-9227}$^{3}$ \\\MakeLowercase{and} Patricia Yanguas \orcidlink{0000-0001-9767-5554}$^{2,\star}$}}
\thanks{$^{\star}$Corresponding author: \texttt{yanguas@unavarra.es}}

\address{\small
$^{1}$Embry-Riddle Aeronautical University, Aerospace Engineering,\\
Daytona Beach, FL 32114, USA\\
$^{2}$Institute for Advanced Materials and Mathematics (INAMAT$^2$),
Department of Statistics, Computer Science and Mathematics,
Public University of Navarre (UPNA),\\ Campus de Arrosadía,
31006 Pamplona, Spain\\
$^{3}$Facultad de Ciencias Jurídicas, Sociales y de la Educación,
Universidad Viña del Mar, Diego Portales 90, Viña del Mar, Chile
}

	

\begin{document}

\begin{abstract}
The 14 Herculis system hosts two massive giant planets on eccentric and mutually inclined orbits, making it a compelling test case for long-term stability in the spatial planetary three-body problem. We develop an analytical description of its dynamics using averaging, symplectic reduction of symmetries and KAM theory. For parameters compatible with the observed architecture, the secular quadrupole Hamiltonian possesses two elliptic relative equilibria that organise the nearby dynamics. The octupole contribution introduces a second slow angle and yields a properly-degenerate multiscale Hamiltonian. Under the corresponding smallness, non-resonance, and non-degeneracy conditions, families of quasi-periodic orbits coming from
the relative equilibria of the quadrupolar approximation fill
invariant $2$-tori. Restoring the two fast mean anomalies and the cyclic angle conjugate to the angular-momentum, and considering the full system, produce two Cantor families of five-dimensional quasi-periodic tori in the reconstructed planetary problem. These invariant structures provide a Hamiltonian mechanism for the long-term confinement of eccentric and mutually inclined configurations compatible with 14 Her.
\end{abstract}

\maketitle

\section{Introduction}

The architectures of multi-planet systems retain valuable evidence of their formation, dynamical evolution, and long-term stability~\cite{Winn2015}. Although many compact systems are characterised by nearly coplanar orbits, a growing number of observed architectures exhibit substantial eccentricities, mutual inclinations, and signatures of strong dynamical excitation~\cite{He2020}. Such dynamically hot systems provide natural laboratories for secular theory, Hamiltonian perturbation methods, and stability mechanisms in the planetary $N$-body problem, while posing a fundamental question: how can eccentric and non-coplanar planetary configurations remain stable over secular timescales?

\begin{figure}
    \centering
    \begin{tikzpicture}[scale=1.3]
	\node[inner sep=0pt] (caso1) at (0,-0.1) 
    {\includegraphics[scale=0.2]{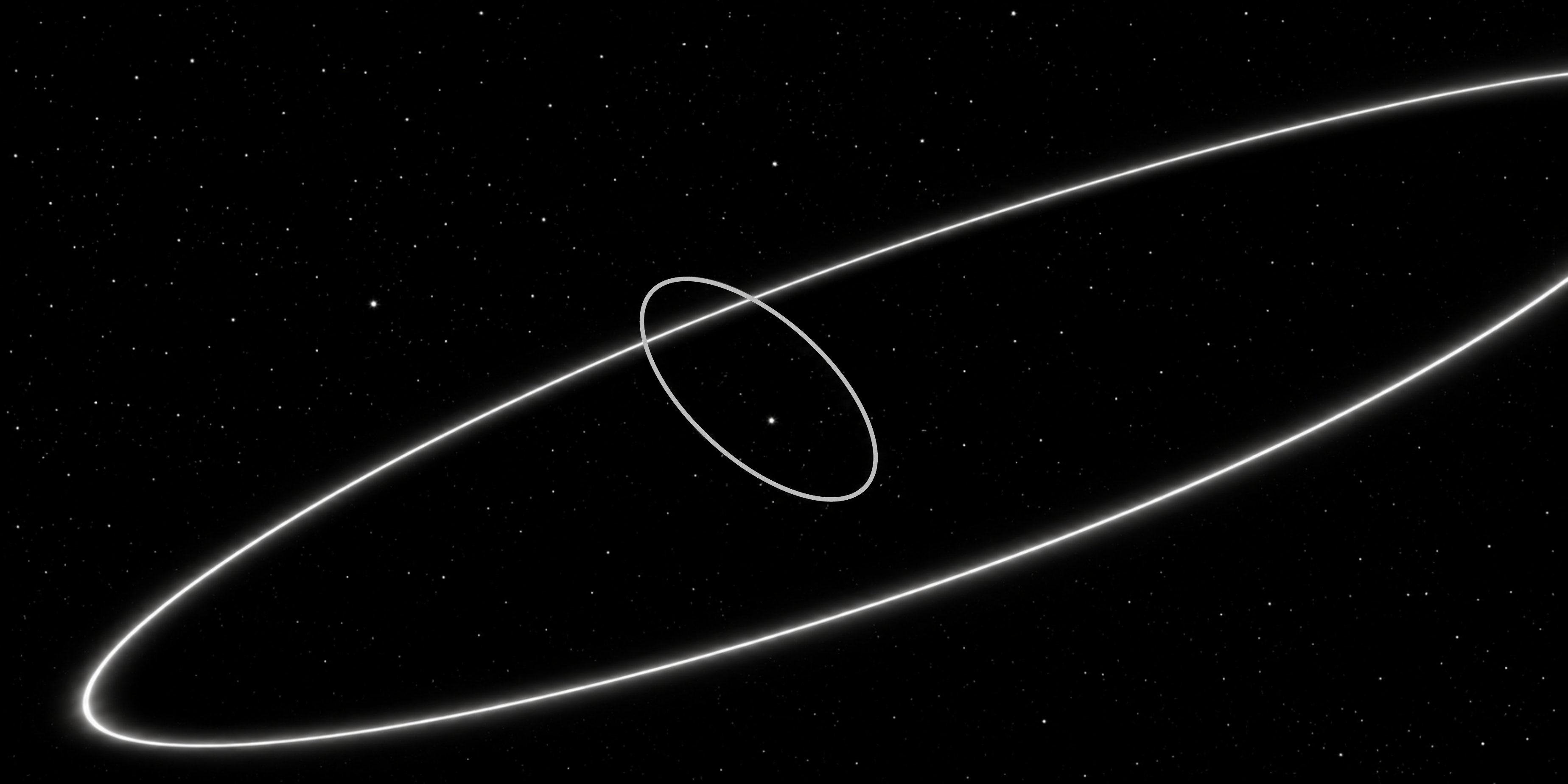}};
    %{\includegraphics[scale=0.135]{Images/Her142.png}};
    %{\includegraphics[scale=0.305]{Images/Her14.pdf}};
	\coordinate [label=below:\textcolor{white}{\small 14 Her}] (E2) at (0,0.3);
    \coordinate [label=below:\textcolor{white}{\small 14 Her b}] (E2) at (-1,1);
    \coordinate [label=below:\textcolor{white}{\small 14 Her c}] (E2) at (1,1.6);
    \end{tikzpicture}
    \caption{Artistic representation of the 14 Herculis exoplanetary system.}
    \label{fig:14Her}
\end{figure}

The nearby K0V star 14 Herculis is a particularly compelling example. It hosts two massive giant planets, 14 Her b and 14 Her c, whose orbital architecture has been progressively constrained through long-baseline radial-velocity measurements, absolute astrometry, and direct imaging. Joint radial-velocity and astrometric analyses established that both planets move on eccentric and mutually inclined orbits, strongly disfavouring coplanarity~\cite{BardalezGagliuffi2021}. JWST/NIRCam coronographic observations subsequently provided a direct image of the cold giant planet 14 Her c, adding a new observational constraint on this dynamically excited system~\cite{BardalezGagliuffi2025}. By incorporating additional Gaia astrometry, \cite{XiaoFeng2025} further constrained the mutual inclination to 35.3$^{+6.8^\circ}_{-7.3}$ (see Fig.~\ref{fig:14Her}). 
Earlier work had already examined dynamically stable two-planet solutions for 14 Her using the radial-velocity constraints then available~\cite{Gozdziewski2008}. More recently, direct $N$-body integrations of current orbital posteriors have assessed the long-term stability and present-day secular evolution of the system, while supporting a primordial scattering scenario involving the ejection of an additional massive planet~\cite{Lu2026}.
Taken together, these results point to a dynamically active past, plausibly involving planet--planet scattering followed by secular angular-momentum exchange and long-term orbital excitation. They also make 14 Her an unusually valuable system in which analytical dynamics can be confronted with a well-constrained three-dimensional architecture.

From a dynamical perspective, 14 Her is naturally modelled as a hierarchical spatial three-body system consisting of a central star, an inner giant planet, and a more distant outer giant planet. The small mass of the planets compared to that of the star makes the two Keplerian motions the leading approximation, with the mutual planetary interaction acting as a perturbation. Averaging over the fast mean anomalies yields a secular Hamiltonian governing the slow evolution of the eccentricities, mutual inclination, and apsidal angles~\cite{Naoz2013}. Although the classical quadrupole-level Kozai--Lidov mechanism provides a useful starting point~\cite{Kozai1962,Lidov1962}, the architecture of 14 Her requires its generalisation to hierarchical systems with finite planetary masses and an eccentric outer orbit, including both quadrupolar and octupolar contributions~\cite{Naoz2013}.
In this setting, higher-order secular expansions and corrections to double averaging provide important benchmarks for assessing truncation accuracy and the possible accumulation of short-period effects over long timescales~\cite{HamersPortegiesZwart2016,LuoKatzDong2016}.

The observed architecture lies beyond the scope of planar and test-particle approximations: both planets have non-negligible masses, their orbital planes are mutually inclined, and the outer orbit is highly eccentric. A consistent long-term description must therefore retain the spatial Hamiltonian coupling among eccentricities, inclinations, and apsidal variables. To reduce rotational symmetry without losing the symplectic structure, we use a coordinate map adapted to the total angular momentum. The two types of variables involved have distinct roles. A Cartesian--nodal set followed by polar--symplectic variables---represented below by the modified polar--nodal, or double-Whittaker, stage---provides an intermediate geometric description of the nodal frame~\cite{Whittaker1904,Deprit1981}. Applying planar Delaunay transformations to the polar canonical pairs then completes Deprit's action--angle variables~\cite{Deprit1983,FerrerOsacar1994a,FerrerOsacar1994b}. In this complete Deprit set, the Delaunay-type pairs isolate the mean anomalies and arguments of periapsis, while the collective nodal actions encode the total angular momentum and its vertical projection; rotational invariance makes their conjugate angles cyclic. Deprit's action--angle variables therefore perform the Jacobi elimination of the nodes symplectically on an open subset of the full phase space and expose the fast mean anomalies required for secular averaging~\cite{PalacianSayasYanguas2013}. A complementary symmetry-reduction framework for the general planetary \(N\)-body problem, formulated in terms of angular-momentum and Laplace--Runge--Lenz invariants, was recently developed by~\cite{ConferenceHer14}. Here we specialise this geometric perspective to the spatial three-body architecture of 14 Herculis and extend it to secular relative equilibria and KAM (Kolmogorov-Arnold-Moser) stability.

The aim of this paper is to identify an analytical mechanism for the long-term stability of configurations compatible with the observed architecture of 14 Her. Rather than refining the orbital fit, we use the measured masses and orbital elements as input for a spatial planetary three-body model. After symmetry reduction and averaging over the two fast angles, we identify two elliptic relative equilibria of the quadrupole Hamiltonian in a region close to the observationally inferred configuration. These equilibria act as the leading-order organising centres of the secular dynamics, with nearby trajectories displaying bounded variations in eccentricity, mutual inclination, and the apsidal variables.

The octupole contribution introduces a second secular angle and breaks the integrability of the quadrupole approximation. Near the equilibria considered here, however, it remains a higher-order contribution to a properly degenerate multiscale Hamiltonian. Under the required smallness, frequency non-resonance, and non-degeneracy hypotheses, KAM theory yields two Cantor families of invariant tori surrounding these configurations. Reinstating the two fast mean anomalies and the cyclic angle associated with angular-momentum reduction then gives five-dimensional quasi-periodic tori in the reconstructed planetary problem. Motions on these tori remain confined for all time, thereby providing a Hamiltonian stability mechanism for eccentric and mutually inclined architectures compatible with the observed behaviour of 14 Her.

This Hamiltonian interpretation complements observational orbit fitting and direct numerical integration. 
A closely related normal-form approach was applied to the $\upsilon$ Andromed{\ae} system, where a librational KAM $2$-torus of the spatial secular model was constructed and its existence established through a computer-assisted proof~\cite{Locatelli2022}; see also the KAM-type analysis performed to some exoplanetary systems in the setting of the spatial three-body problem~\cite{Caracciolo2024}.

Our objective is not to reconstruct one particular invariant torus with high numerical precision, but to identify the geometric structures that organise a family of stable motions within the region of phase space allowed by the observations. This perspective is particularly useful for exoplanetary systems, whose orbital parameters remain subject to significant uncertainties. The method also provides a route for studying other hierarchical two-planet systems with eccentric and mutually inclined orbits, such as HAT-P-11, which contains a close-in super-Neptune and an eccentric outer super-Jupiter on an orbit significantly inclined relative to that of the inner planet~\cite{AnEtAl2025}. It also suggests extensions to systems containing more than two planets, such as Kepler-129, a hierarchical three-planet system comprising two close-in sub-Neptunes and a long-period giant planet inferred to be inclined relative to the inner planetary plane~\cite{ZhangEtAl2021}.

The paper is organised as follows. Section~\ref{sec:hamiltonian_model} introduces the spatial planetary Hamiltonian and its associated secular system. Section~\ref{subsec:quadrupole_equilibria} identifies the quadrupole relative equilibrium associated with the 14 Her architecture and analyses its local dynamics. Section~\ref{subsec:hanliyi} addresses the persistence and reconstruction of the corresponding invariant tori and their associated quasi-periodic motions, while numerical illustrations are presented in Section~\ref{sec:numerical}. 
Section~\ref{sec:conc} summarises the results and discusses their scope. In Appendix~\ref{Appen}, we explain how the angle coordinates must be modified in some of the transformations we perform to ensure that they are symplectic.

%Opción 1
\begin{figure*}
    \centering\large
    \begin{tikzpicture}[scale=1]
	\node[inner sep=0pt] (caso1) at (0,0) {\includegraphics[scale=0.6]{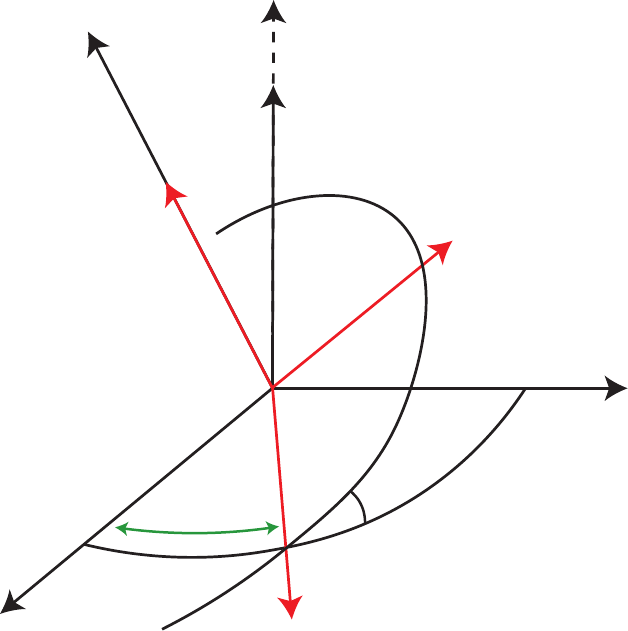}};
	\coordinate [label=below:\textcolor{black}{$\ell$}] (E2) at (0,-2.5);
    \coordinate [label=below:\textcolor{black}{\small $I$}] (E2) at (0.3,-1.82);
    \coordinate [label=below:\textcolor{black}{$\nu$}] (E2) at (-1.2,-1.8);
    \coordinate [label=below:\textcolor{black}{$j$}] (E2) at (3,0);
    \coordinate [label=below:\textcolor{black}{$i$}] (E2) at (-3,-3);
    \coordinate [label=below:\textcolor{black}{$B$}] (E2) at (-0.75,3.4);
    \coordinate [label=below:\textcolor{black}{$k$}] (E2) at (-0.75,2.4);
    \coordinate [label=below:\textcolor{black}{$m$}] (E2) at (1.6,1.2);
    \coordinate [label=below:\textcolor{black}{$G$}] (E2) at (-2.6,3);
    \coordinate [label=below:\textcolor{black}{$n$}] (E2) at (-1.8,1.6);
    \end{tikzpicture}
    \qquad
    \begin{tikzpicture}[scale=1]
	\node[inner sep=0pt] (caso1) at (0,0) {\includegraphics[scale=0.8]{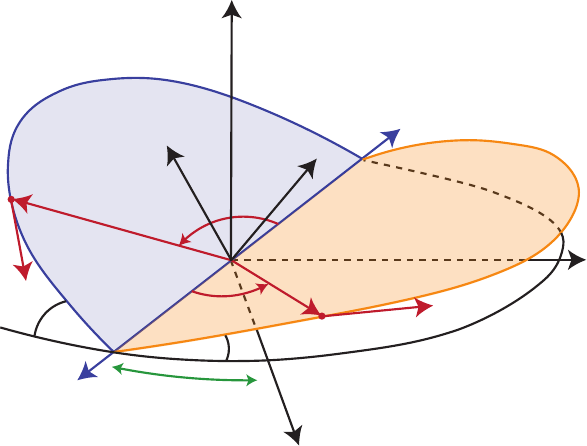}};
	\coordinate [label=below:\textcolor{black}{$\ell$}] (E2) at (0.3,-2.5);
    \coordinate [label=below:\textcolor{black}{$\nu_1$}] (E2) at (-1.5,-2.1);
    %\coordinate [label=below:\textcolor{black}{$\nu_2$}] (E2) at (2.5,-1.6);
    \coordinate [label=below:\textcolor{black}{\normalsize $I_1$}] (E2) at (-1.1,-1.47);
    \coordinate [label=below:\textcolor{black}{\normalsize $I_2$}] (E2) at (-3.1,-1.1);
    \coordinate [label=below:\textcolor{black}{\normalsize $\ell_1$}] (E2) at (-3.1,-2);
    \coordinate [label=below:\textcolor{black}{\normalsize $q_2$}] (E2) at (-3.6,0.8);
    \coordinate [label=below:\textcolor{black}{\normalsize $p_2$}] (E2) at (-3.9,-0.2);
    \coordinate [label=below:\textcolor{black}{\normalsize $n_1$}] (E2) at (-1.7,1.4);
    \coordinate [label=below:\textcolor{black}{\normalsize $n_2$}] (E2) at (0.5,1.2);
    \coordinate [label=below:\textcolor{black}{\normalsize $n$}] (E2) at (-0.8,3.4);
    \coordinate [label=below:\textcolor{black}{\normalsize $\ell_2=-\ell_1$}] (E2) at (1.5,1.7);
    \coordinate [label=below:\textcolor{black}{\normalsize $\theta_1$}] (E2) at (-0.7,-0.5);
     \coordinate [label=below:\textcolor{black}{\normalsize $\theta_2$}] (E2) at (-0.6,0.5);
     \coordinate [label=below:\textcolor{black}{\normalsize $q_1$}] (E2) at (0.5,-1.3);
    \coordinate [label=below:\textcolor{black}{\normalsize $p_1$}] (E2) at (2.1,-0.9);
    \coordinate [label=below:\textcolor{black}{\normalsize $m$}] (E2) at (4.2,-0.3);
    \end{tikzpicture}
    \caption{Left: inertial frame
\(\{\mathbf{i},\mathbf{j},\mathbf{k}\}\) and invariable frame \(\{\pmb{\ell},\mathbf{m},\mathbf{n}\}\). Right: geometry of the angular momentum vectors and the Deprit variables. }
    \label{fig:DiagramaG}
\end{figure*}

\section{The spatial planetary three-body problem}
\label{sec:hamiltonian_model}

%%%%%%%%%%%%%%%%

\subsection{Heliocentric formulation}

We consider the exoplanetary system 14 Herculis as a spatial planetary three-body problem consisting of a dominant central star of mass $m_0$ and two planets of masses $\varepsilon m_1$ and $\varepsilon m_2$, where $\varepsilon \ll 1$ measures the planetary character of the system. Let $(\mathbf{x}_i,\mathbf{y}_i)\in\mathbb{R}^3\times\mathbb{R}^3$ denote the (inertial) barycentric position--momentum pairs of the star ($i=0$) and the two planets, with $i=1$ and $i=2$ corresponding to the inner and outer planets, respectively. Let $(\mathbf{q}_i,\mathbf{p}_i)\in\mathbb{R}^3\times\mathbb{R}^3$, $i=1,2$, denote the related heliocentric position--momentum pairs of the planets. 
After fixing the centre of mass at the star and reducing the translational symmetry, the configuration variables are defined by $\mathbf{q}_i = \mathbf{x}_i-\mathbf{x}_0$, $i=1,2$. Thus, the origin of the reduced configuration space is located at the centre of the star, while the conjugate variables $\mathbf{p}_i$ are the associated planetary 
momenta~\cite{Meyer1999,MeyerOffin2017,ChierchiaPinzari2011}. To keep the transformation symplectic, the new momenta have to be scaled with respect to the inertial ones, so that $\mathbf{p}_i=\mathbf{y}_i/\varepsilon$ for $i=1,2$. In the symplectic heliocentric coordinates, the Hamiltonian takes the form

\begin{equation}
H = H_{\mathrm{Kep}} + \varepsilon H_{\mathrm{pert}},
\label{eq:Hfull}
\end{equation}
where
\begin{equation}
\begin{split}
H_{\mathrm{Kep}} =
\sum_{i=1}^{2}
\left(\frac{\|\mathbf{p}_i\|^2}{2\mu_i} - \frac{\mu_i M_i}{\|\mathbf{q}_i\|}\right),
\quad
H_{\mathrm{pert}} =
\frac{\mathbf{p}_1\cdot\mathbf{p}_2}{m_0} - \frac{m_1 m_2}{\|\mathbf{q}_1-\mathbf{q}_2\|},
\end{split}
\label{eq:Hpert}
\end{equation}
with
\[
\mu_i = \frac{m_0 m_i}{m_0 + \varepsilon m_i},
\quad
M_i = m_0 + \varepsilon m_i,
\quad
i = 1,2.
\]

The Keplerian part describes two uncoupled elliptic orbits revolving around the star at the common focus, whereas the perturbation accounts for the mutual interaction between the planets.

For the 14 Herculis system, we adopt the stellar mass as the mass unit, so that $m_0=1$. The posterior median dynamical masses reported by~\cite{BardalezGagliuffi2025} are $0.97\pm0.04\,M_\odot$ for the star, $8.9^{+1.3}_{-1.7}\,M_{\rm J}$ for 14 Her b, and $7.9^{+1.6}_{-1.2}\,M_{\rm J}$ for 14 Her c. Their central values correspond to the dimensionless planet-to-star mass ratios $8.76\times10^{-3}$ and $7.77\times10^{-3}$, respectively. Since both ratios are close to $10^{-2}$, we adopt $\varepsilon=10^{-2}$ as a convenient common mass scale, or bookkeeping parameter, and write them as $\varepsilon m_1$ and $\varepsilon m_2$, with $m_1 \simeq 0.876$ and $m_2 \simeq 0.777$. Thus, $\varepsilon$ isolates the common small mass scale while the rescaled masses $m_1$ and $m_2$ remain of order unity. This scaling makes the planetary regime explicit: the leading star--planet Keplerian dynamics is contained in $H_{\mathrm{Kep}}$, whereas the planet--planet coupling is smaller by a factor $\varepsilon$ and appears through $\varepsilon H_{\mathrm{pert}}$. Such a dominant-central-mass hierarchy is characteristic of systems formed by one star and its planets and provides the starting point for the qualitative analysis developed below.

%\textcolor{red}{Decir de qu\'e tamaño es eps para Herculis 14 y que en realidad es el responsable de tratar el problema en el r\'egimen planetario, que es lo usual para sistemas formados por una estrella y sus planetas. Decirlo con un poco de detalle ya que es nuestro punto de partida}.

The 14 Herculis system is also strongly hierarchical in orbital size. Using the central values of the semi-major axes obtained by~\cite{XiaoFeng2025},
\[
a_1 \simeq 2.84\ {\rm au},
\qquad
a_2 \simeq 28.1\ {\rm au},
\]
the ratio of the inner to the outer semi-major axis is
\[
\delta = \frac{a_1}{a_2} \simeq 0.101 \simeq 0.1.
\]
%\textcolor{red}{El valor de $\delta$ y su relaci\'on con $\varepsilon$ (ponerlo) tambi\'en es clave, ya que nos da la idea de separar el t\'ermino octupolar del cuadrupolar. Hay que decir que $\delta$ est\'a relacionado con el cocente de los semi-ejes mayores de las \'orbitas de los dos planetas alrededor de la estrella. El hecho de que el planeta exterior est\'e alejado de su estrella en comparaci\'on con el interior permite tratar el problema cuadrupolar+octupolar como lo hacemos, facilitando su an\'alisis.}
Thus, the outer planet orbits at a characteristic distance about ten times larger than that of the inner planet. The parameters $\varepsilon$ and $\delta$ describe two distinct physical hierarchies: $\varepsilon$ measures the small planetary masses relative to the star, whereas $\delta$ measures the separation of the two orbital scales. For 14 Herculis, the adopted values satisfy the approximate ordering $\delta\simeq\varepsilon^{1/2}$. Treating $\delta$ as the bookkeeping parameter for the orbital hierarchy, the Newtonian interaction term admits the Legendre expansion
\begin{equation}
\frac{1}{\|\mathbf{q}_1 - \mathbf{q}_2\|} = \frac{1}{r_2} \sum_{n\ge 0} \left(\frac{r_1}{r_2}\right)^n P_n(\cos\Phi),
\label{eq:LegendreExpansion}
\end{equation}
where
\[
r_i = \|\mathbf{q}_i\|,\qquad i = 1,2, \quad
\cos\Phi = \frac{\mathbf{q}_1 \cdot \mathbf{q}_2}{r_1 r_2}.
\]

The expansion converges whenever $r_1 < r_2$, a condition satisfied throughout the region of phase space relevant to the 14 Herculis system. Retaining terms through fourth order in $\delta$ yields the monopole, dipole, quadrupole and octupole contributions. This is equivalent to expanding explicitly (\ref{eq:LegendreExpansion}) up to $n=3$. Adding the dot product of the momenta $\mathbf{p}_i$ (indirect kinetic term), the perturbation can be written as
\[
H_{\rm pert}
= \mathbf{p}_1 \cdot \mathbf{p}_2 + H_{\rm mono} + H_{\rm dipo} + H_{\rm quad} + H_{\rm octu} + {\mathcal O}(\delta^5),
\]
where $H_{\rm mono} = {\mathcal O}(\delta)$, $H_{\rm dipo} = {\mathcal O}(\delta^2)$, $H_{\rm quad} = {\mathcal O}(\delta^3)$ and $H_{\rm octu}={\mathcal O}(\delta^4)$. The indirect kinetic term scales as ${\mathcal O}(\delta^{1/2})$. Since the complete perturbation is of the form $\varepsilon H_{\rm pert}$, the corresponding contributions are ordered as
\[
\varepsilon  \mathbf{p}_1 \cdot \mathbf{p}_2 = {\mathcal O}(\varepsilon\delta^{1/2}) = 
{\mathcal O}(\varepsilon^{5/4}), \]\vspace*{-0.35cm}
\[
\varepsilon H_{\rm mono} = {\mathcal O}(\varepsilon\delta)
= {\mathcal O}(\varepsilon^{3/2}),
\qquad
\varepsilon H_{\rm dipo} = {\mathcal O}(\varepsilon\delta^2)
={\mathcal O}(\varepsilon^2),
\]
and
\[
\varepsilon H_{\rm quad} = {\mathcal O}(\varepsilon\delta^3)
= {\mathcal O}(\varepsilon^{5/2}),
\qquad
\varepsilon H_{\rm octu} = {\mathcal O}(\varepsilon\delta^4)
= {\mathcal O}(\varepsilon^3).
\]
Numerically, the associated bookkeeping scales of the quadrupole and octupole terms are approximately $10^{-5}$ and $10^{-6}$, respectively. The octupole contribution is therefore smaller than the quadrupole contribution by an additional factor of $\delta\simeq0.1$. This separation of scales justifies retaining the quadrupole and octupole terms while treating higher multipoles as smaller corrections (of order ${\mathcal O}(\varepsilon^{7/2})$, i.e. with bookkeeping scale $10^{-7}$), and makes the hierarchical problem amenable to the analysis developed below. Nevertheless, it is important to bear in mind that, when the eccentricity of the outer planet is sufficiently large, the octupole-order terms may become comparable in magnitude to the quadrupole-order terms, and similarly with the higher-order terms; see, for instance, \cite{Naoz2016,Naoz2020}. This will need to be taken into account in the following sections.

\subsection{Reduction of the nodes and Deprit's variables}
\label{subsec:deprit}
Following~\cite{Deprit1983}, the implementation used by~\cite{FerrerOsacar1994a,FerrerOsacar1994b}, and the broader symmetry-reduction framework for planetary systems developed in~\cite{ConferenceHer14}, we perform the reduction of the rotational symmetry by introducing coordinates adapted to the total angular momentum vector. We first pass from the heliocentric Cartesian variables to a suitable set of polar-nodal type, and then to two pairs of action-angle coordinates, called Deprit's variables. 

Let $\mathbf{G}_i = \mathbf{q}_i \times \mathbf{p}_i$, $i = 1,2$ be the angular momentum vectors of the two planets, and the total angular momentum
\[
\mathbf{G} = \mathbf{G}_1 + \mathbf{G}_2,
\qquad
G = \|\mathbf{G}\|,
\qquad
B = \mathbf{G}\cdot\mathbf{k},
\]
where \(\mathbf{k}\) is the fixed vertical unit vector of the inertial frame $\{\mathbf{i}, \mathbf{j}, \mathbf{k}\}$.
The modulus \(G\) and the projection \(B\) are constants of motion associated with the
rotational symmetry.

The reduction of the nodes is performed in several stages. First, we move from the inertial reference frame to the invariable frame $\{\pmb{\ell},\mathbf{m}, \mathbf{n}\}$ defined by $\mathbf{G}$, see Figure~\ref{fig:DiagramaG} left. More precisely, when $\mathbf{G} \neq \mathbf{0}$ is not parallel to $\mathbf{k}$, we define the orthonormal basis
\[ \mathbf{n} = \dfrac{\mathbf{G}}{G}, \quad \pmb{\ell} = \dfrac{\mathbf{k} \times \mathbf{n}}{\|\mathbf{k} \times \mathbf{n} \|},\quad \mathbf{m} = \mathbf{n} \times \pmb{\ell}.\]

Next we introduce two other frames. Assume $\mathbf{G}_i \neq \mathbf{0}$ for $i = 1, 2$. There is a unique polar decomposition $\mathbf{G}_i = G_i \mathbf{n}_i$ with $\| \mathbf{n}_i \| = 1$. We define the angle $I_i$ such that $\mathbf{n} \cdot \mathbf{n}_i = \cos I_i$ with $0 \le I_i \le \pi$. If $I_i \in ( 0, \pi )$ there exists a unique direction $\pmb{\ell}_i$ with $\mathbf{n} \times \mathbf{n}_i = \pmb{\ell}_i \sin \, I_i$, $\| \pmb{\ell}_i \| = 1$. Specifically $I_i$ is the angle between $\mathbf{G}$ and $\mathbf{G}_i$, that is, the inclination of the ellipse $i$ with respect to the invariable plane. 

The nodal direction $\pmb{\ell}_i$ lies in the invariable plane $\{\pmb{\ell}, \mathbf{m}\}$ and can be decomposed as $\pmb{\ell}_i = \pmb{\ell} \cos \nu_i + \mathbf{m} \sin \, \nu_i$ with $0 \le \nu_i < 2 \pi$ the nodal angles. By construction, $\pmb{\ell}_2=-\pmb{\ell}_1$, and hence $\nu_2 = \nu_1 + \pi$ mod $2\pi$. For $i = 1,2$, the nodal frames are defined through the three orthonormal directions $\{ \mathbf{n}_i, \pmb{\ell}_i, \mathbf{m}_i \}$, where $\mathbf{m}_i = \mathbf{n}_i \times \pmb{\ell}_i$. 

The inertial and invariable frames are related by means of the rotation matrices ${\textsf R}_1 = {\mathcal R}_3(\nu) {\mathcal R}_1(I)$ ${\mathcal R}_3(\nu_1) {\mathcal R}_1(I_1)$, 
${\textsf R}_2 = {\mathcal R}_3(\nu) {\mathcal R}_1(I)$ ${\mathcal R}_3(\nu_1 + \pi) {\mathcal R}_1(I_2)$, where $I$ is the inclination of the invariable plane orthogonal to $\mathbf{G}$ determined by $\cos I = B/G$, and $\nu$ is the angle between $\mathbf{i}$ and $\pmb{\ell}$ called argument of the node, and is determined by $\cos \nu = \mathbf{i} \cdot \pmb{\ell}$ and $\sin \nu = \mathbf{j} \cdot \pmb{\ell}$. 

Now the first step of the transformation is
\[ (\mathbf{q}_1, \mathbf{q}_2, \mathbf{p}_1, \mathbf{p}_2) \rightarrow (\mathbf{Q}_1, \mathbf{Q}_2, \mathbf{P}_1, \mathbf{P}_2), \]
given by $\mathbf{q}_i = {\textsf R}_i \cdot \mathbf{Q}_i$ and $\mathbf{p}_i = {\textsf R}_i \cdot \mathbf{P}_i$ for $i = 1, 2$, with the feature that the third components of $\mathbf{Q}_i$, $\mathbf{P}_i$ are zero. The next step consists of a modified double Whittaker transformation \cite{Whittaker1904}, also named symplectic polar variables \cite{Deprit1981} applied on each planar pair $(\mathbf{Q}_i,\mathbf{P}_i)$:
\[ (\mathbf{Q}_1, \mathbf{Q}_2, \mathbf{P}_1, \mathbf{P}_2) \rightarrow 
(r_1, \theta_1, \nu_1, r_2, \theta_2, \nu, R_1, \Theta_1, G, R_2, \Theta_2, B), \]
where $r_i$ is the radial distance of the body describing the ellipse $i$ while $R_i$ corresponds to its conjugate momentum, $\theta_i$ is the so-called argument of the latitude with conjugate variable $\Theta_i = G_i$. The set $\{ r_i$, $\theta_i$, $\nu_1$, $\nu$, $R_i$, $\Theta_i$, $G$, $B \}$ corresponds to Deprit's polar-nodal variables. Note that $G$ and $B$ are respectively conjugate to $\nu_1$ and $\nu$.

At this point, for each orbital plane we introduce planar Delaunay variables
\[
\{ l_i, g_i, L_i, G_i \}, \qquad i = 1,2,
\]
where \(l_i\) is the mean anomaly, \(g_i\) the argument of periapsis, \(L_i\) is the action conjugate to \(l_i\), and \(G_i\) is the modulus of the angular momentum $\mathbf{G}_i$. The angles \(g_1\) and \(g_2\) play the role of arguments of periapsis in the orbital planes of the two ellipses. The angle $g_i$ is measured from the mutual node $\pmb{\ell}_i$. The passage from $\{ r_i, \theta_i, R_i, \Theta_i \}$ is the standard one. Thus the full coordinate system becomes
\[
\{ l_1, g_1, \nu_1, l_2, g_2, \nu, L_1, G_1, G, L_2, G_2, B \}.
\]
The coordinates $l_i$, $L_i$, $G_i$ coincide with the spatial Delaunay coordinates whereas $g_i$, $\nu_1$, $\nu$, $G$ and $B$
are specific to the coordinate system introduced by \cite{Deprit1983}. The whole transformation remains symplectic.

Considering the semi-major axis $a_i$ of each planet, and its corresponding eccentricity $e_i$, the variables satisfy
\[
L_i = \mu_i \sqrt{M_i a_i},
\qquad
G_i = L_i \sqrt{1 - e_i^2},
\qquad i = 1,2.
\]

The inclinations of the orbital planes with respect to the invariable plane can be obtained by means of
\[ \cos I_1 = \frac{G^2 + G_1^2 - G_2^2}{2 G G_1}, \qquad \cos I_2 = \frac{G^2 + G_2^2 - G_1^2}{2 G G_2}, \]
whereas $\sin I_i$ is determined as $\sqrt{1-\cos^2 I_i}$.

The mutual inclination between the two orbital planes is determined from the previous formulae, and is given by
\begin{equation}
\cos (I_1+I_2) =
\frac{G^2-G_1^2-G_2^2}{2G_1G_2}.
\label{eq:mutualinclination}
\end{equation}
%%%%%%%%

In these coordinates the Keplerian part of the Hamiltonian depends only on the actions $L_i$:
\begin{equation}
H_{\mathrm{Kep}}
=
-\frac{M_1^2\mu_1^3}{2L_1^2}
-\frac{M_2^2\mu_2^3}{2L_2^2},
\label{eq:HKepDeprit}
\end{equation}
while the perturbation depends on all variables but the nodes:
\[
H_{\mathrm{pert}}
=
H_{\mathrm{pert}}
(l_1, l_2, g_1, g_2, L_1, L_2, G_1, G_2; G).
\]
Thus, Hamiltonian \eqref{eq:Hfull} becomes independent of the nodal angles $\nu$ and $\nu_1$ and independent of $B$ in this set of variables. Therefore, $\nu$, $B$ and $G$ are first integrals and the system is reduced from six to four degrees of freedom.

%%%%%%%%%%%%%%%%%
\iffalse
\begin{figure}
\centering\large
    \begin{tikzpicture}[scale=1]
	\node[inner sep=0pt] (caso1) at (0,0) {\includegraphics[scale=0.8]{Images/triangulo.pdf}};
	\coordinate [label=below:\textcolor{black}{$G_1$}] (E2) at (-1.5,-2.3);
    \coordinate [label=below:\textcolor{black}{$G_2$}] (E2) at (0.8,2.8);
    \coordinate [label=below:\textcolor{black}{$G$}] (E2) at (3.8,0.25);
    \coordinate [label=below:\textcolor{black}{$e$}] (E2) at (-1,0.35);
    \coordinate [label=below:\textcolor{black}{$e$}] (E2) at (1,0);
    \coordinate [label=below:\textcolor{black}{$f$}] (E2) at (-2,-0.1);
    \coordinate [label=below:\textcolor{black}{$f$}] (E2) at (2.2,0.55);
    \coordinate [label=below:\textcolor{black}{$I_2$}] (E2) at (-2.3,0.45);
    \coordinate [label=below:\textcolor{black}{$I_1$}] (E2) at (-2.8,-0.15);
    \end{tikzpicture}
    \captionsetup{width=\linewidth}
    \caption{Caption}
    \label{Figura1}
\end{figure}
\fi
%%%%%%%%%%%%%%%%%%

\subsection{Averaging of the fast angles}
\label{subsec:averaging}

The Hamiltonian \eqref{eq:Hfull} expressed in Deprit's action-angle variables still contains the fast Keplerian motions associated with the mean anomalies \(l_1\) and \(l_2\). Since the 14 Herculis system is sufficiently far from low-order mean-motion resonances, a first-order Lie normalisation can be performed with respect to these angles. The associated canonical transformation is $\varepsilon$-close to the identity and removes the dependence on \(l_1\) and \(l_2\) up to first order in the perturbation. The hierarchical architecture of the system prevents the appearance of small denominators associated with commensurabilities of the mean motions, thereby ensuring the validity of the normalisation procedure in a neighbourhood of the observed configuration. This averaging procedure can be understood as a first-order canonical transformation within the framework of Lie transformations.

After dropping terms higher than the octupolar ones, the secular Hamiltonian is therefore obtained as
\begin{equation}
\bar{H} = 
\frac{1}{4\pi^2}
\int_0^{2\pi}\!\!\!
\int_0^{2\pi}
H_{\mathrm{pert}} \, dl_1 \, dl_2.
\label{eq:doubleaverage}
\end{equation}

In order to compute the average explicitly, it is convenient to use different
anomalies for the two integrations. For the inner orbit we employ the eccentric
anomaly \(E_1\), since the perturbing terms contain powers of the radial
distance \(r_1\) in the numerator. The averaging is performed using the classical expression for the differential of the mean anomaly in terms of the eccentric and true 
anomalies. Using the partial derivative of $E_1$ with respect to $l_1$
and that of $f_2$ with respect to $l_2$ given by
\begin{equation}
\frac{\partial E_1}{\partial l_1} = \frac{1}{1-e_1\cos E_1},\quad \frac{\partial f_2}{\partial l_2} = \frac{a_2^2\eta_2}{r_2^2} = \frac{(1+e_2 \cos f_2)^2}{\eta_2^3},
\label{eq:l1average}
\end{equation}
where $\eta_2 = \sqrt{1-e_2^2}$, the double-averaged Hamiltonian is therefore computed through the successive quadratures
\[
\bar H = \frac{\eta_2^3}{4\pi^2} \int_0^{2\pi} \!\!\! \int_0^{2\pi}
H_{\mathrm{pert}} \, (1-e_1\cos E_1) (1 + e_2 \cos f_2)^{-2} \, d E_1 d f_2.
\]
In the process, the argument of latitude is written as $\theta_i = f_i + g_i$. Moreover, standard formulae are applied to relate the momentum $R_1$ to $\sin E_1$ and $r_1$ to $\cos E_1$, while expressing $\cos f_1$ and $\sin f_1$ in terms of $\cos E_1$ and $\sin E_1$ and, finally, $r_2$ as a function of $\cos f_2$. The final expression inside the integrals is a finite Fourier series in $E_1$, $f_2$, thus the integrals are straightforward to calculate.

After averaging and truncating the transformation's tail, the mean anomalies disappear from the Hamiltonian and the actions $L_i$ become constant. The resulting Hamiltonian is the so-called secular Hamiltonian up to the octupolar approximation (once the remainder is dropped) and depends only on the slow angles and their conjugate momenta, namely
\begin{equation}
\begin{aligned}
H_{\rm sec}
&=
H_{\rm Kep}
+
\varepsilon
\left(
\bar{H}_{\rm mono}(L_2)
+
\bar{H}_{\rm quad}(g_1, G_1, G_2; L_1, L_2, G) \right.
\\ & \left. \hspace*{1.53cm} + \,
\bar{H}_{\rm octu}(g_1, g_2, G_1, G_2; L_1, L_2, G)
\right).
\end{aligned}
\label{eq:Hsec}
\end{equation}

The averages of the terms $\mathbf{p}_1 \cdot \mathbf{p}_2$ and $H_{\rm dipo}$ with respect to $l_1$ and $l_2$ are zero.

The quadrupole contribution defines an integrable one-degree-of-freedom subsystem, whereas the octupole term introduces the second secular angle \(g_2\), producing a genuine two-degree-of-freedom Hamiltonian system. Truncation at higher orders also produces a Hamiltonian system
of two degrees of freedom.

%%%%%%%%%%%%%%%%%%%
\subsection{Monopole, quadrupole and octupole secular Hamiltonians}
\label{subsec:quad_oct}

The monopole contribution obtained from the double-averaging procedure described in the previous subsection is not relevant for the analysis of the secular Hamiltonian in the reduced space, but we shall need it when applying the KAM theorem in Section~\ref{subsec:hanliyi}. The double average gives
\begin{equation}
\bar{H}_{\rm mono} = \kappa, \quad \mbox{with} \quad
\kappa = -\frac{m_1 m_2 M_2 \mu_2^2}{L_2^2}.
\label{eq:Hmono}
\end{equation}
The quadrupole contribution is given by
\begin{equation}
\bar{H}_{\rm quad} = \tilde{\kappa}\, h_{\rm quad},
\label{eq:Hquad}
\end{equation}
where $\tilde{\kappa}$ is a positive constant depending only on the masses and the semi-major axes
\[
\tilde{\kappa} =
\frac{m_1 m_2 M_2^{3} \mu_2^{6} L_1^{2}}
     {64 M_1^{2} \mu_1^{4} L_2^{3}}
\]
and $h_{\rm quad}$ is given by
\begin{equation}
\begin{split}
h_{\rm quad}
=
\frac{1}{G_1^{2} G_2^{5}}
\Big[
& -3(G^{2} - G_1^{2})^{2}(5L_1^{2} - 3G_1^{2}) \\
& +G_{2}^2 (6 G^2 - 2 G_{1}^2 - 3 G_{2}^2) (5L_1^{2} - 3G_1^{2}) \\
& -15\left((G^2 + G_{1}^2 + G_{2}^2)^2 - 2(G^4 + G_{1}^4 + G_{2}^4) \right)
\\
&\times(L_1^2 - G_1^2) \cos(2g_1) \Big].
\end{split}
\label{eq:hquad}
\end{equation}
As pointed out above, a fundamental feature of \eqref{eq:Hquad} is that it depends only on the pair \( (g_1, G_1) \), while \(L_1\), \(L_2\), \(G_2\) and \(G\) remain constant parameters. Consequently, the quadrupole approximation defines an integrable one-degree-of-freedom Hamiltonian on the reduced space ${\mathcal T}_{L_1, L_2, G_2}$ introduced in \cite{PalacianSayasYanguas2013}.

The next contribution in the hierarchy is the octupole term,
which is written as
\begin{equation}
\bar{H}_{\rm octu} = \kappa'\, h_{\rm octu},
\label{eq:Hoct}
\end{equation}
with
\[
\kappa' =
-\frac{15 m_1 m_2 M_2^4 \mu_2^8 L_1^3}{2048 M_1^3 \mu_1^6 L_2^4}
\]
and
\[
\begin{aligned}
h_{\rm octu} = {}&\frac{\sqrt{(L_1^2 - G_1^2)(L_2^2 - G_2^2)}}{G_1^3 G_2^8}
\Big\{
70G_1 G_2 {\mathcal A} \cos^3g_1 \cos g_2 \\
& +G_1 G_2 \left[{\mathcal B} + 105 {\mathcal A} \cos(2g_1) \right] \cos g_1 \cos g_2 \\
& -105 {\mathcal C} (L_1^2 - G_1^2) \cos^2g_1 \sin g_1 \sin g_2 \\
& +\frac{1}{2}(G^2 - G_1^2 - G_2^2)
\left[{\mathcal D} + 35 {\mathcal A} \cos(2g_1)\right] \sin g_1 \sin g_2
\Big\},
\end{aligned}
\]
where
\[
\begin{aligned}
{\mathcal A} = {}& -(L_1^2 - G_1^2) \left(G^4 - 2G^2(G_1^2 + G_2^2) + (G_1^2 - G_2^2)^2 \right),\\
{\mathcal B} = {}& -135G_1^6 + 174G_1^4 G_2^2 - 135G_1^2 G_2^4
+175L_1^2 G_1^4 - 126L_1^2G_1^2 G_2^2 +175L_1^2 G_2^4 + 5G^4(35L_1^2 - 27G_1^2)
\\& -10G^2(G_1^2 + G_2^2)(35L_1^2 - 27G_1^2), \\
{\mathcal C} = {}& G^6 - 3G^4(G_1^2 + G_2^2)
-(G_1^2 - G_2^2)^2(G_1^2 + G_2^2)  +G^2(3G_1^4 + 2G_1^2 G_2^2 + 3G_2^4),\\
{\mathcal D} = {}& -125G_1^6 + 154G_1^4 G_2^2 - 125G_1^2 G_2^4
+245L_1^2 G_1^4 - 266L_1^2 G_1^2 G_2^2 +245L_1^2 G_2^4 + 5G^4(49L_1^2 - 25G_1^2)
\\
& -10G^2(G_1^2 + G_2^2)(49L_1^2 - 25G_1^2).
\end{aligned}
\]
The relative orders are $\bar{H}_{\rm mono} = {\mathcal O}(\delta)$, $\bar{H}_{\rm quad} = {\mathcal O}(\delta^3)$ and
$\bar{H}_{\rm octu} = {\mathcal O}(\delta^4)$.

Unlike the quadrupole contribution,
\(\bar{H}_{\rm octu}\) depends explicitly on the angle
\(g_2\), thereby producing a genuine two-degree-of-freedom
Hamiltonian system. Hence, the octupolar approximation of
the secular Hamiltonian is defined on a reduced space 
whose dimension is four. In \cite{PalacianSayasYanguas2013} it is named
${\mathcal S}_{L_1, L_2, G}$. Its geometry and dynamics are analysed in \cite{PalacianSayasYanguas2015}. 
More precisely it corresponds to a symplectic manifold with singularities for
certain combinations of the parameters $L_1$, $L_2$ and $G$.

A direct inspection of the averaged octupole term (\ref{eq:Hoct}) shows that the
dependence on \(g_2\) appears through trigonometric factors
proportional to \(\sin g_2\) and \(\cos g_2\). 

\section{Two elliptic relative equilibria of the quadrupole secular Hamiltonian}
\label{subsec:quadrupole_equilibria}

\subsection{Relative equilibria of the quadrupole model}

The relative equilibria of the double-averaged quadrupole Hamiltonian and their linear stability have a well-established history. 
Early treatments based on Delaunay variables include the construction of near-circular quasi-periodic motions by~\cite{JefferysMoser1966} and Harrington's averaged model of the hierarchical stellar three-body problem~\cite{Harrington1969}.
Within the same framework, \cite{LidovZiglin1976} gave a systematic account of the relative equilibria, their stability, and the associated bifurcation curves.
\cite{FerrerOsacar1994a,FerrerOsacar1994b} subsequently revisited Harrington's Hamiltonian using Deprit's variables instead of Delaunay elements, thereby placing the elimination of the nodes in a symplectic coordinate framework and refining the analysis of the parameter plane. 
This approach was extended by~\cite{PalacianSayasYanguas2013}, who combined regular and singular reduction to incorporate the singular strata of the fully reduced phase space and to clarify the occurrence of relative equilibria, linear stability, and bifurcations associated with them.
\cite{PalacianSayasYanguas2015} then reconstructed the corresponding invariant KAM tori associated to the elliptic relative equilibria. These tori are of dimension five for the Hamiltonian system associated to (\ref{eq:Hfull}). Building on these results, the purpose of this section is not to repeat the general classification, but to identify the elliptic branch relevant to the observed architecture of 14 Herculis and to use it as the organising centre for the local and KAM analyses developed below.

The quadrupole Hamiltonian \eqref{eq:Hquad} defines an integrable
one-degree-of-freedom system on the reduced phase space
\(\mathcal{T}_{L_1, L_2, G_2}\). Since the quantities \(L_1\), \(L_2\),
\(G_2\), and \(G\) are constant, the dynamics is governed by the
canonical pair \((g_1, G_1)\), whose equations of motion are
\begin{equation}
\dot g_1 =
\frac{\partial \bar{H}_{\rm quad}}{\partial G_1},
\qquad
\dot G_1 =
-\frac{\partial \bar{H}_{\rm quad}}{\partial g_1}.
\label{eq:quadrupole_equations}
\end{equation}

The relative equilibria of the quadrupole approximation correspond the
critical points of \(\bar{H}_{\rm quad}\) on ${\mathcal T}_{L_1, L_2, G_2}$. The analysis, including the cases in which either of the two trajectories is circular, 
the two orbits are coplanar, or the inner motion is rectilinear, was carried out 
in~\cite{PalacianSayasYanguas2013}. Here, we focus on the motion of the 14 Herculis 
system and therefore restrict our analysis to
\[
0 < G_i < L_i,
\qquad
|G-G_2| < G_1 < G+G_2.
\]

Since the dependence of \(\bar{H}_{\rm quad}\) on the angular variable \(g_1\) appears only through the term \(\cos(2g_1)\), the second equation in \eqref{eq:quadrupole_equations} yields $\sin(2g_1) = 0$, and therefore
\begin{equation}
g_1 \in \left\{0, \,{\pi}/{2}, \,\pi, \,{3\pi}/{2} \right\}.
\label{eq:g1candidates}
\end{equation}
Substituting \eqref{eq:g1candidates} into the first equation of
\eqref{eq:quadrupole_equations} produces an algebraic equation for the
angular momentum \(G_1\). The branches $g_1 = 0$, and $g_1 = \pi$ lead to one family of equilibria, whereas $g_1 = {\pi}/{2}$, and $g_1 = {3\pi}/{2}$ generate a second family, which will be of particular interest in the present work.

For $g_1 \in \{ \pi/2, 3\pi/2 \}$ the equilibrium condition reduces to
\begin{equation}
5L_1^2 (G_2^4 + G^4) - G_1^4 (5L_1^2 - 8G_1^2 + 4G_2^2)
-2G^2 (5L_1^2 G_2^2 + 4G_1^4)
= 0.
\label{eq:G1equilibrium}
\end{equation}
Equation \eqref{eq:G1equilibrium} is a cubic polynomial in \(G_1^2\). It corresponds to the polynomial
$s_6(G_1)$ in \cite{PalacianSayasYanguas2013}, p. 172. Its physically admissible roots determine the quadrupolar relative equilibria. 
The relative equilibria obtained from \eqref{eq:G1equilibrium} act as organising centres for the surrounding secular dynamics. In the next subsection we evaluate this equation using the observational parameters of the 14 Herculis system and identify the equilibrium configuration that
best matches the observed planetary architecture.

\subsection{The 14 Herculis equilibrium}
\label{subsec:14her_equilibrium}

We now specialise the quadrupole model to the observed 14 Herculis system. Throughout this work we adopt the orbital parameters reported in the recent analyses based on radial velocity, astrometric, and direct-imaging observations \cite{BardalezGagliuffi2021,BardalezGagliuffi2025,XiaoFeng2025}. The values used in our computations are
\begin{equation}
\begin{aligned}
a_1 &= 2.84\ {\rm au},
&
a_2 &= 28.1\ {\rm au},
\\
e_1 &= 0.371,
&
e_2 &= 0.64,
\\
\varepsilon m_1 &= 8.9\,M_{\rm J},
&
\varepsilon m_2 &= 7.9\,M_{\rm J},
\\
i_1 &= 147.3^\circ,
&
\\
\Omega_1 &= 197.8^\circ,
&
\Omega_2 &= 205.1^\circ,
\end{aligned}
\label{eq:14HerData}
\end{equation}
where \(\Omega_1\) and \(\Omega_2\) denote the longitudes of the ascending nodes, two of the Delaunay angles. The angle $i_1$ represents the inclination of the inner ellipse with respect to the plane spanned by $\{ \mathbf{i}$, $\mathbf{j} \}$. The inclination \(i_2\) of the outer planet is left as a free parameter to be determined from the equilibrium condition.

At this stage, we emphasise that the values given above correspond to the parameters of the true orbits. In what follows, however, we will work with their averaged values, which introduces a discrepancy of order ${\mathcal O}(\varepsilon)$. This does not affect our main achievements. 

The orbital parameters \eqref{eq:14HerData} determine the actions and angular momenta
\[
L_i = \mu_i \sqrt{M_i a_i}, \quad G_i = L_i \sqrt{1 - e_i^2},
\qquad i = 1,2.
\]
Therefore all coefficients appearing in the quadrupole Hamiltonian \eqref{eq:Hquad} are fixed except for the dependence on the total angular momentum \(G\). To express \(G\) in terms of the orbital elements, we introduce the projections of the angular momentum vectors onto the inertial vertical axis, 
\[
H_i = G_i \cos i_i,
\qquad i = 1,2,
\]
that correspond to the Delaunay actions conjugate to the nodes $\Omega_i$. Using the geometry of the angular momentum triangle, the modulus of the total angular momentum is given by
\begin{equation}
G^2 = G_1^2 + G_2^2 + 2H_1 H_2 + 2\sqrt{G_1^2 - H_1^2} \sqrt{G_2^2 - H_2^2} \cos(\Omega_1 - \Omega_2).
\label{eq:Gtotal}
\end{equation}
Since the inclination \(i_1\) of the inner planet and the nodal longitudes \(\Omega_1\) and \(\Omega_2\) are quantities known from the observational fit, equation \eqref{eq:Gtotal} expresses the total angular momentum \(G\) as a function of the unknown inclination \(i_2\). Consequently, the equilibrium condition \eqref{eq:G1equilibrium} may be regarded as an equation involving the unknown inclination \(i_2\).

Substituting the observational parameters \eqref{eq:14HerData} together with the expression \eqref{eq:Gtotal} into \eqref{eq:G1equilibrium} yields the inclination
\[
i_2^0 = 98.6^\circ.
\]
This value differs from the inclination \(i_2 = 114^\circ\) reported in the observational analyses \cite{BardalezGagliuffi2025} or \( i_2 = {111.9^{+5.4}_{-4.5}}^\circ \) in \cite{XiaoFeng2025}. However, the purpose of the present construction is not to reproduce exactly the observed orbital configuration, but rather to identify a relative equilibrium of the quadrupole Hamiltonian lying in its neighbourhood. The difference of approximately \(13^\circ\) between the inclination of Xiao and Feng and ours remains small when compared with the characteristic angular scales of the secular dynamics and places the observed system within the region of phase space organised by the equilibrium. Consequently, the equilibrium provides a dynamically meaningful reference configuration that is fully compatible with the current observational constraints on the 14 Herculis system.
Additionally, the difference \(i_1-i_2^0\) yields a mutual inclination of \(48.7^\circ\), which differs from the reported 
value of \({35.3^{+6.8}_{-7.3}}^\circ\) in \cite{XiaoFeng2025}. Applying formula (\ref{eq:mutualinclination}) gives a mutual inclination around \(49.05^\circ\), which is compatible to \(48.7^\circ\). The corresponding value of the total angular momentum is \(C = 4.17783\), which remains constant throughout the analysis.

It is important to emphasise that the observed planetary system is not expected to lie exactly at the equilibrium point. Rather, the equilibrium acts as an organising centre of the secular dynamics. Trajectories initiated in its neighbourhood remain confined to a small region of the reduced phase space and exhibit oscillations of their orbital elements around the equilibrium values. In this sense, the equilibrium provides a natural reference configuration around which the long-term dynamics of the system should be studied.

This observation motivates the analysis of the local dynamics in a neighbourhood of this equilibrium, which will be the starting point for the construction of the quasi-periodic solutions associated with the 14 Herculis system.

%\textcolor{red}{Los p\'arrafos anteriores quedan bien. Tiene que quedar claro que nuestro punto de vista es partir de un equilibrio relativo en el sistema secular formado por el hamiltoniano cuadrupolar, de forma que dicho equilibrio da lugar a soluciones casi-peri\'odicas del sistema de partida. Como la parte octupolar es importante, tenemos que ``ajustar" dicho equilibrio para analizar la din\'amica del sistema cuadrupolar+octupolar y probar que el equilibrio relativo da lugar a toros de dimensi\'on dos en el espacio reducido \({\mathcal S}_{L_1, L_2, G}\) del cual se deducen los toros de dimensi\'on 5 en el problema de partida. Finalmente las soluciones casi-peri\'odicas de dichos toros explican la din\'amica del sistema Herculis.}

\subsection{Local dynamics and multiscale structure}
\label{subsec:local_dynamics}

The two equilibrium points identified in the previous subsection provides a
natural reference configuration for the secular dynamics of the
14 Herculis system. In order to study motions in their
neighbourhood, we introduce local coordinates centred at the
equilibria,
\begin{equation}
x = g_1 - \frac{\pi}{2} \,\,\, \mbox{or},\,\,\, x = g_1 - \frac{3\pi}{2}
\qquad
y = G_1 - G_1^0,
\label{eq:localcoordinates}
\end{equation}
where \( G_1^0 = 1.36469 \) is the solution of the equilibrium equation
\eqref{eq:G1equilibrium} associated with the observed
architecture.

Since ${\partial \bar{H}_{\rm quad}}/{\partial g_1} = {\partial \bar{H}_{\rm quad}}/{\partial G_1} = 0$ at the equilibrium, the Taylor expansion of the quadrupole Hamiltonian takes the form
\begin{equation}
\bar{H}_{\rm quad} = \tilde{\kappa} \bar{H}_{\rm quad}^0 + \tilde{\kappa} \left( a\, x^2 + b\, y^2 \right) + {\mathcal O}\!\left( \|(x,y)\|^3 \right),
\label{eq:quadraticexpansion}
\end{equation}
where $\bar{H}_{\rm quad}^0$ and the coefficients \(a\) and \(b\) depend on the parameters \(L_1\), \(L_2\), \(G_2\), \(G\) and \(G_1^0\). Moreover, $\bar{H}_{\rm quad}^0$ is independent of $x$ and $y$. We get:
\[
\begin{aligned}
\bar{H}_{\rm quad}^0 &= \frac{2}{G_s G_2^5}
\left[
12G_s
\left(
(G^2 - G_s)^2
-(2G^2 + G_s) G_2^2 + G_2^4
\right) \right. \\ 
&  \left. \qquad -5 L_1^2
\left(
3(G^2 - G_s)^2
-2(3G^2 + G_s) G_2^2
+3G_2^4
\right)
\right], \\
a &= \frac{30}{G_s G_2^5} (L_1^2 - G_s) \left(G^4 + (G_s - G_2^2)^2 - 2G^2 (G_s + G_2^2) \right), \\
b &= -\frac{6}{G_s^2 G_2^5}
\left(
-24G_s^3
+15L_1^2 G^4
+15L_1^2 G_2^4
+G_s^2 (5L_1^2 + 4G_2^2) \right. \\
& \left. \qquad +2G^2(4G_s^2 - 15L_1^2 G_2^2)
\right),
\end{aligned}\]
where $G_s = (G_1^0)^2$. 
The numerical values corresponding to the 14 Herculis system show that the equilibrium is elliptic if $a b > 0$ and this happens when $G_1^0$ lies in the interval $(1.02351, 1.46957)$, as well as other intervals that are not relevant to our analysis (considering the actions $L_1$, $G$, and $G_2$ fixed). Consequently, the linearised dynamics consists of small oscillations around the equilibrium point and the nearby trajectories are organised by closed invariant curves of the quadrupole approximation.

\begin{figure}
\centering
\captionsetup{width=\linewidth}
\includegraphics[width=0.7\textwidth]{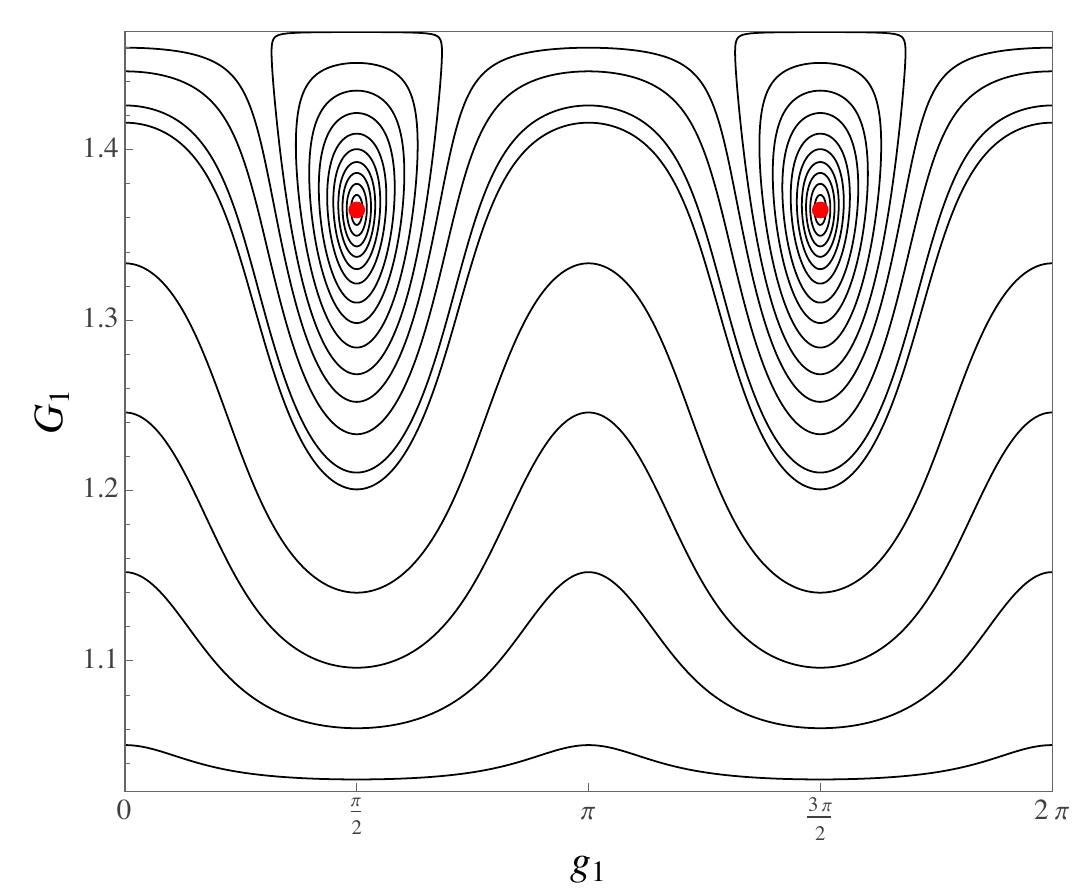}
\caption{Phase portrait of the quadrupole Hamiltonian in the projection of the reduced phase space
\((g_1, G_1)\). The black curves are level sets of the quadrupole Hamiltonian,
while the blue curve corresponds to the energy level passing through the
relative equilibrium. The red point marks the equilibrium.}
\label{fig:PhasePortraitQuad}
\end{figure}

Figure~\ref{fig:PhasePortraitQuad} illustrates the phase portrait of the quadrupole approximation in the neighbourhood of the equilibrium determined in Section~\ref{subsec:14her_equilibrium}.

The secular Hamiltonian, up to terms of order \(\delta^4\), is obtained after incorporating the octupole contribution,
\[
\bar{H}_{\rm qo} = \bar{H}_{\rm quad} + \bar{H}_{\rm octu}.
\]
As discussed in Section~\ref{subsec:quad_oct}, the octupole terms are one order higher in the hierarchy parameter \(\delta\) and contain the additional secular angle \(g_2\). However, when the eccentricity of the outer body is significant, as in the 14 Herculis system, the octupolar terms can be of the same order as the quadrupolar ones, as already pointed out in several papers, see for instance \cite{Naoz2020} and references therein.

Indeed, for the two-degree-of-freedom Hamiltonian system defined by the secular 
terms up to octupolar order, elliptic relative equilibria (i.e., linearly stable 
equilibria) would give rise to invariant $5$-tori and quasi-periodic motions in the full system. The existence of such equilibria as functions of the system parameters was analysed in \cite{Naoz2020}. However, the 14 Herculis system does not fall within the parameter range considered in that study for their numerical approach. Consequently, the desired KAM tori cannot be obtained from these equilibria. 

Alternatively our approach consists in establishing the existence of invariant 2-tori in the 
reduced space ${\mathcal S}_{L_1, L_2, G}$, where the secular Hamiltonian $\bar{H}_{\rm qo}$ is defined, as invariant structures surrounding a periodic orbit of the system. 
This periodic orbit is obtained by perturbing the equilibrium $g_1 = \pi/2$ or 
$3\pi/2$, $G_1 = G_1^0$ through the inclusion of the octupolar terms. Consequently, $g_1(t)$ is of the form $\pi/2 + {\mathcal O}(\delta)$ or $3\pi/2 + {\mathcal O}(\delta)$, $g_2 = g_2(t) \in [0, 2\pi)$, $G_2 = G_2(t)$, and $G_1 = G_1(t)$, with $G_1(t)$ representing a $\delta$-perturbation of the equilibrium value $G_1^0$.

To construct the local normal form, we restrict our analysis to a region of the reduced 
space ${\mathcal S}_{L_1, L_2, G}$ in which the quadrupolar terms dominate, albeit slightly, over the octupolar terms. At the same time, we require $g_1$ to remain close to $\pi/2$
or $3\pi/2$, and $G_1$ to remain close to its equilibrium value $G_1^0$.

This separation of scales allows the secular Hamiltonian to be regarded as a multiscale perturbation of the elliptic equilibrium. The quadrupole Hamiltonian provides the leading integrable approximation, while the octupole contribution appears at a higher order in the hierarchy parameter and introduces a second slow angular variable. As a consequence, the local dynamics near the equilibrium possesses the hierarchical structure required for the application of KAM theory.

In the next section we exploit this multiscale Hamiltonian structure to establish the persistence of invariant tori and to construct quasi-periodic motions associated with the 14 Herculis system.

\section{Quasi-periodic motions and KAM tori}
\label{QuasiKAM}

\subsection{A KAM theorem suitable for Hamiltonian systems with high-order degeneracy}
\label{subsec:HLY}

To establish the persistence of quasi-periodic motions we employ a theorem due to \cite{HanLiYi2010}, which is a KAM-type result specifically designed for Hamiltonian systems with several asymptotic scales and proper degeneracies. In its simplest form, the theorem states that if a multiscale Hamiltonian can be written in action-angle coordinates $(J, \psi)$ as
\begin{equation}
H(J,\varphi) = h_0(J) + \delta^{\alpha_1} h_1(J) + \cdots + \delta^{\alpha_r} h_r(J) + P(J, \psi),
\label{eq:HanLiYiForm}
\end{equation}
with $\alpha_i$ positive rational numbers arranged in increasing order,  \(P\), the perturbation, sufficiently small, smaller than the integrable part depending only 
on $J$ and the corresponding non-degeneracy conditions, that will be tackled explicitly below, are satisfied then, a Cantor family of invariant tori persists under perturbation. The Han–Li–Yi theorem extends Arnold’s classical approach to the persistence of KAM tori in properly degenerate systems, where the unperturbed Hamiltonian is independent of some of the action variables, while the perturbation removes this degeneracy. Arnold’s result applies to systems in which the unperturbed Hamiltonian is formulated as a two-scale Hamiltonian system. 

The Han–Li–Yi theorem also generalises earlier results of Bruno and R\"ussmann, in particular R\"ussmann’s approach, which requires the weakest assumptions and asserts that the frequency vector associated with the integrable part of the Hamiltonian should not lie in any $(n-1)$-dimensional subspace. Indeed, when a KAM theorem is applied by taking into account only two different scales, the perturbation must be very small in order to draw conclusions about the persistence of solutions, as has been done in many applications to the three- and $N$-body problems.

The Han--Li--Yi theorem has been applied in the context of the spatial three-body problem. 
In particular, \cite{PalacianSayasYanguas2015} considers elliptic relative equilibria in different reduced spaces associated with the quadrupolar approximation and uses them to reconstruct invariant $5$-tori and the corresponding quasi-periodic motions of the problem in several configurations. Among these, some of the equilibria are studied in the reduced space ${\mathcal T}_{L_1,L_2,G_2}$, including those considered in the present work; see, in particular, case (a) of Table 6.1 in \cite{PalacianSayasYanguas2015}. In that setting, however, the octupolar contribution is assumed to be of higher order than the quadrupolar one and therefore does not affect the conclusions. In the present work, we require additional conditions on the Hamiltonian system to establish the persistence of the KAM tori arising from the quadrupolar approximation. These conditions are studied in the following subsections.

\subsection{Arranging the Hamiltonian in a suitable portion of phase space}
\label{subsec:arrangeh}
After the symplectic reductions and transformations performed above, it becomes clear that the Hamiltonian obtained for the 14 Herculis system fits naturally within this framework. Nevertheless, when considering the associated secular Hamiltonian, the situation is more subtle. Although the quadrupole approximation yields an integrable Hamiltonian possessing the elliptic equilibria identified in Section~\ref{subsec:14her_equilibrium}, the octupole contribution enters at higher order in $\delta$. However, it can become comparable in magnitude to the quadrupole contribution, even within a small neighbourhood of the two equilibrium points.

To ensure that the main contribution of the octupole terms remain of higher order than the dominant quadrupole terms when $g_1$ is close to $\pi/2$ or $3\pi/2$, and $G_1$ is close to $G_1^0$, we make use of some freedom we have in the choice of $G_2$. Specifically, we slightly shift $G_2$ away from its observational value, obtained from (\ref{eq:14HerData}), namely $3.15432$. It is sufficient to consider
$G_2 = 3.15432 + \Delta_{G_2}$, with $\Delta_{G_2}$ varying between $1/16$ and $1/12$. For this choice of $G_2$, the corresponding value of $G_1$ obtained by solving (\ref{eq:G1equilibrium}), hence using that $g_1 = \pi/2$, is approximately in the interval $[1.3092, 1.3233]$. In the estimation we have employed the values of $L_1$, $G$ deduced from (\ref{eq:14HerData}). The analysis with $g_1$ in a neighbourhood of $3\pi/2$ is essentially the same. We then introduce the local coordinates
\begin{equation}
x = g_1 - \frac{\pi}{2},
\qquad
y = G_1 - G_1^0,
\label{eq:localcoordinates2}
\end{equation}
where, in contrast to the previous section, $G_1^0$ is no longer a constant, but rather the expression obtained from the appropriate solution of (\ref{eq:G1equilibrium}). Thus, $G_1^0$ depends not only on $G_2$, but also on $L_1$ and $G$. Viewed as a transformation in the phase space associated with the Hamiltonian \eqref{eq:Hfull}, the change of coordinates \eqref{eq:localcoordinates2} is canonical provided that the angles $g_2$, $l_1$, and $\nu_1$ are modified accordingly. The explicit corrections to these angles are given in Appendix \ref{Appen}. We denote the transformed angles by $\bar{g}_2$, $\bar{l}_1$ and $\bar{\nu}_1$. The remaining coordinates, namely $L_1$, $l_2$, $L_2$, $G_2$, $G$, $\nu$, and $B$, do not change.

Next we introduce a pair of action-angle coordinates by means of
\begin{equation}
    x = \sqrt{2 J/\omega} \sin \varphi, \qquad y = \sqrt{2 \omega J} \cos \varphi,
\end{equation}
where
\begin{equation} 
\omega = 
\frac{G_1^0 \sqrt{5 (L_1^2 - G_s)
\left((G + G_2)^2 - G_s\right) \left(G_s - (G - G_2)^2 \right)}}{
\sqrt{5L_1^2 \left(3(G^2 - G_2^2)^2 + G_s^2\right) + 4(2G^2 + G_2^2 - 6G_s) G_s^2
}
},
\label{freq}
\end{equation}
and $G_s = (G_1^0)^2$. Note that $\omega=\sqrt{a/b}$ where $a$, $b$ are introduced in Section~\ref{subsec:local_dynamics} but now they depend on the action $G_2$. When $L_1$ and $G$ take values close to the ones corresponding to the 14 Herculis equilibrium, while $G_2$ remains near $3.15432 + \Delta_{G_2}$, with $\Delta_{G_2}$ in 
$[1/16, 1/12]$, the square roots appearing in $\omega$ are positive. The value of $\omega$ lies approximately in $[0.2734, 0.2892]$. These bounds on \(\omega\) hold uniformly throughout the portion of phase space considered in our analysis.

To guarantee that the transformation remains entirely symplectic, we need to modify the angles $\bar{l}_1$, $\bar{\nu}_1$, and $\bar{g}_2$, denoting the resulting variables by 
$\tilde{l}_1$, $\tilde{\nu}_1$, and $\tilde{g}_2$, respectively. The necessary corrections are presented in Appendix \ref{Appen}. All other variables, that is, 
$L_1$, $l_2$, $L_2$, $G_2$, $G$, $\nu$, and $B$, remain unchanged.

Applying the transformations just introduced and performing a Taylor expansion around $x=y=0$, the perturbation becomes

\[ 
\begin{aligned}  
\bar{H}_{\rm quad}^a =\,\,& \frac{m_1 m_2 M_2^3 \mu_2^6 L_1^2}{32 M_1^2 \mu_1^4 L_2^3 G_2^5 G_s} \, {\mathcal Q}_a, \\
\bar{H}_{\rm quad}^b =\,\,& -\frac{3\sqrt{5} m_1 m_2 M_2^3 \mu_2^6 L_1^2 \sqrt{L_1^2 - G_s}}{16 M_1^2 \mu_1^4 L_2^3 G_2^5 G_s^{3/2}} \, {\mathcal Q}_b,
\end{aligned}
\]
with 
\[
\begin{aligned}  
{\mathcal Q}_a =\, & 
12G_s
\left(
(G^2 - G_s)^2
-(2G^2 + G_s) G_2^2 + G_2^4
\right)  \\ 
&  \quad - 5 L_1^2
\left(
3(G^2 - G_s)^2
-2(3G^2 + G_s) G_2^2
+3G_2^4
\right), \\
{\mathcal Q}_b =\, & \sqrt{2 G^2 (G_s + G_2^2) - G^4 - (G_s - G_2^2)^2 } \\ & \quad \times \sqrt{5L_1^2 \left(3(G^2 - G_2^2)^2 + G_s^2\right) + 4(2G^2 + G_2^2 - 6G_s) G_s^2
} \, J. 
\end{aligned}\]
The square roots appearing in the definition of ${\mathcal Q}_b$ remain positive throughout the domain of the actions $L_1$, $G$, and $G_2$ under consideration. The subsequent terms in $\bar{H}_{\rm quad}$ are multiplied by monomials $J^{j/2}$, $j\in \mathbb{N}$, starting at $J^{3/2}$, and can be regarded as a Fourier series in $\varphi$. For sufficiently small $J$, these terms remain very small compared with ${\mathcal Q}_b$ and the Fourier series is convergent.

Concerning $\bar{H}_{\rm octu}$ it has a term independent of $J$ given by
\[ \bar{H}_{\rm octu}^a = - \frac{15 m_1 m_2 M_2^4 \mu_2^8}{512 M_1^3 \mu_1^6 L_2^4 G_2^8 G_s^{3/2}} \, {\mathcal O}_a, \]
where 
\[ 
\begin{aligned}  
{\mathcal O}_a =\, & {\mathcal G} \left[ 4G_s \left(5 (G^2 - G_s)^2 - (10 G^2 + 7 G_s) G_2^2 + 5 G_2^4 \right) \right. \\ & \left. \, - 
 7 L_1^2 \left(5 (G^2 - G_s)^2 - 2 (5 G^2 + 3 G_s) G_2^2 + 5 G_2^4 \right) \right] \sin g_2, \\
\mbox{with } {\mathcal G} =\, & L_1^3 (G_s + G_2^2 - G^2) \sqrt{(L_1^2 - G_s) (L_2^2 - G_2^2)}.
 \end{aligned}
 \]
The rest of terms in $\bar{H}_{\rm octu}$ are expressible as a double Fourier series (finite in $g_2$) and $\varphi$ whose coefficients are multiplied by powers of $J^{j/2}$, $j \in \mathbb{N}$ with $j \ge 1$, thus the series is convergent is $J$
is chosen small enough. Notice that $g_2$ should be still written in terms of $\tilde{g}_2$. We call $\bar{H}_{\rm octu}^b$ those terms in $\bar{H}_{\rm octu}$ multiplied by exponents of $J^{1/2}$ and
$J$, but avoid writing it as it is a lengthy expression.

As
\[
J = \frac{(G_1 - G_1^0)^2 + \frac{1}{4}(2g_1 - \pi)^2\omega^2}{2\omega},
\]
the quantity \(J\) measures the squared distance of the point \((g_1, G_1)\) from \((\pi/2, G_1^0)\), up to the scaling determined by \(\omega\). Thus, \(J\) can be made arbitrarily small by taking \(g_1\) sufficiently close to \(\pi/2\) and \(G_1\) sufficiently close to \(G_1^0\). Choosing $0 < J < 0.0182$ implies that
$G_1$ is near $G_1^0$ while $g_1$ is close to \(\pi/2\). More specifically, from the inequalities $|G_1 - G_1^0| \le \sqrt{2\omega J}$, $|g_1 -\pi/2| \le \sqrt{2J/\omega}$, we find 
\[ |G_1 - G_1^0| \le 0.0997 \quad \mbox{and} \quad |g_1 -\pi/2| \le 0.3547. \]

Taking the masses values those of the 14 Herculis equilibrium and the actions $L_1$, $L_2$ also those of the 14 Herculis equilibrium, the values of $\bar{H}_{\rm quad}^a$ are approximately in $[0.000039, 0.000047]$ while the absolute values of $\bar{H}_{\rm quad}^b$ are upper-bounded by $0.000015$, thus $3.1  |\bar{H}_{\rm quad}^b| \approx \bar{H}_{\rm quad}^a$. Thus we estimate $|\bar{H}_{\rm quad}^b|/\bar{H}_{\rm quad}^a \approx \delta^{1/2}$. The terms multiplied by $J^{j/2}$ with $j \ge 3$ are of order $10^{-6}$ or lower, but always of smaller magnitude than $\bar{H}_{\rm quad}^b$.

Regarding the comparison of $\bar{H}_{\rm octu} = \bar{H}_{\rm octu}^a + \bar{H}_{\rm octu}^b$ with $\bar{H}_{\rm quad}^a$, we note that $\bar{H}_{\rm octu}^a$ depends on $\tilde{g}_2$ and $\varphi$, with both angles taking values in $[0, 2\pi)$. Its magnitude, $|\bar{H}_{\rm octu}^a|$, oscillates between zero and a maximum value ranging from $0.0000273$ to $0.0000284$, while its average values lie approximately between $0.000017$ and $0.000018$. Consequently, the ratio $|\bar{H}_{\rm octu}^a|/\bar{H}_{\rm quad}^a$ can be considered of order $\delta^{2/5}$. By contrast, the terms in $\bar{H}_{\rm octu}^b$ can reach magnitudes of up to $10^{-6}$ as $J$ approaches $0.0182$, but remain smaller in magnitude than $\bar{H}_{\rm octu}^a$.

At this point, we can average $\bar{H}_{\rm octu}^a$ with respect to $\tilde{g}_2$, taking $\bar{H}_{\rm quad}^a$ as the unperturbed term. Indeed, $\bar{H}_{\rm octu}^a$ depends on the angles $\tilde{g}_2$ and $\varphi$, whereas $\bar{H}_{\rm quad}^a$ is independent of $J$. Thus, the corresponding transformation is considerably simple. Since the average of 
$\bar{H}_{\rm octu}^a$ with respect to $\tilde{g}_2$ vanishes --all terms in $\bar{H}_{\rm octu}^a$ are of the form $\sin(\tilde{g}_2 + F)$ with $F$ a function that could depend on $G_2$, $\varphi$ and $J$ but not on $\tilde{g}_2$-- the associated homological equation is
\[ \frac{\partial \bar{H}_{\rm quad}^a}{\partial G_2} \frac{\partial W}{\partial \tilde{g}_2} = \bar{H}_{\rm octu}^a,\]
where $W$ denotes the generating function associated with the Lie transformation. This function is readily determined by first computing the antiderivative of $\bar{H}_{\rm octu}^a$ with respect to $\tilde{g}_2$ (which is simply a finite Fourier series in this angle) and then dividing the resulting expression by the partial derivative of $\bar{H}_{\rm quad}^a$ with respect to $G_2$. The latter quantity does not vanish in the region of phase space considered in this analysis. The values that the derivative takes in the portion of phase space that we are considering is always above $0.0003$. 

As a consequence of the transformations performed above, within the region ${\mathcal S}_{L_1, L_2, G}$ determined previously in terms of $G_2$, $g_1$, and $G_1$, the perturbation
$\bar{H}_{\rm quad} + \bar{H}_{\rm octu}$ can be decomposed into the main part $\bar{H}_{\rm quad}^a$, its $\delta^{1/2}$-order perturbation $\bar{H}_{\rm quad}^b$, and the remaining terms arising from $\bar{H}_{\rm quad}$ and $\bar{H}_{\rm octu}$. The latter terms may depend on $\tilde{g}_2$ and $\varphi$, but they should be considered as a perturbation of $\bar{H}_{\rm quad}^a$ entering only at order ${\mathcal O}(\delta^{4/5})$ or higher.

We can say that the ellipticity of the equilibria for the quadrupole approximation guarantees the existence of local action-angle variables, while the hierarchy of scales established in Section~\ref{subsec:local_dynamics} provides the required separation between the dominant and perturbative dynamics.

We now extend the transformations constructed for the secular Hamiltonian to the full Hamiltonian with six degrees of freedom. Consequently, the actions $L_i$ are no longer fixed but are allowed to vary slightly while the angles $l_i$ arise at order ${\mathcal O}(\varepsilon^3)$. It remains to verify the rank condition required by the Han–Li–Yi theorem.

Denote by \(h_0, \ldots, h_4\) the successive
integrable contributions of the multiscale normal form. We have arrived at
\[
\begin{aligned}  
h_0(L_1) =\, &  -\frac{M_1^2 \mu_1^3}{2L_1^2}, \\
h_1(L_1, L_2) =\, & -\frac{M_2^2 \mu_2^3}{2L_2^2}, \\
h_2(L_1, L_2) =\, & -\frac{m_1 m_2 M_2 \mu_2^2}{L_2^2}, \\
h_3(L_1, L_2, G, G_2) =\,\, & \bar{H}_{\rm quad}^a, \\
h_4(L_1, L_2, G, G_2, J) =\,\, & \bar{H}_{\rm quad}^b, 
\end{aligned}
\]
such that in a suitable portion of the phase space $\mathbb{R}^3\times\mathbb{R}^3$ where $H$ introduced in (\ref{eq:Hfull}) reads
\begin{equation} 
\label{HamNor}
H = h_0 + h_1 + \varepsilon h_2 + \varepsilon h_3 + \varepsilon h_4 + {\mathcal O}(\delta^{29/5}), 
\end{equation}
where $h_1$, $h_2$ are of order ${\mathcal O}(\delta)$ while $h_3$ is of order ${\mathcal O}(\delta^3)$ and finally, $h_4$ is
of order ${\mathcal O}(\delta^{7/2})$, all with respect to $h_0$. 

Notice that the Hamiltonian $H$ in normal form is consistent with the discussion of the different orders of the perturbative expansion presented at the end of Section~\ref{sec:hamiltonian_model}. The double averaging over the mean anomalies, the introduction of different sets of canonical coordinates, the restriction to a small neighbourhood of phase space around the two relative equilibria of the quadrupolar terms, and the subsequent averaging over the angle $\tilde{g}_2$ push the non-essential terms to order ${\mathcal O}(\delta^{29/5})$, while retaining the relevant contributions at the preceding orders. It is worth noting that the dependence of $H$ on the action $J$ first appears at order ${\mathcal O}(\varepsilon \delta^{7/2})$, which is slightly below ${\mathcal O}(\delta^{29/5})$. This contribution originates from the quadrupolar perturbation. 

The multiscale structure of $H$ in \eqref{HamNor} becomes apparent at this stage, with four distinct scales emerging and giving rise to invariant tori characterised by slow and fast frequencies across the different scales.

\subsection{Application of the Han-Li-Yi theorem}
\label{subsec:hanliyi}
The Hamiltonian is now in a form suitable for applying the theorem of Han, Li, and Yi. It remains to prove that, starting from the frequencies obtained 
by differentiating the integrable part of $H$ with respect to the five actions, we can construct a matrix of rank five, equal to the dimension of the invariant 
tori we seek to establish.

Order the action variables as
\[
\mathcal{I} = (L_1, L_2, G, G_2, J).
\]
The frequency map is built as the vector
\[ \Omega(\mathcal{I}) = \left( 
\frac{\partial h_0}{\partial L_1} \,,\, 
\frac{\partial h_1}{\partial L_2} \,,\,
\frac{\partial h_2}{\partial L_2} \,,\, 
\frac{\partial h_3}{\partial G} \,,\, 
\frac{\partial h_3}{\partial G_2} \,,\, 
\frac{\partial h_4}{\partial J}\right). \]
The corresponding action blocks have dimensions $n_i - n_{i-1}$, $i = 0, \ldots, 4$, with 
\[
n_0 = 1, \qquad n_1 = 2, \qquad n_2 = 2, \qquad
n_3 = 4, \qquad n_4 = 5.
\]

A convenient non-zero $5\times5$ minor of the matrix formed by the block-frequency map and all its first-order derivatives is
\begin{equation}
\mathcal{R}_{5} =
\begin{pmatrix}
\frac{\partial h_0}{\partial_{L_1}} & 0 & 0 & 0 & 0\\
0 & \frac{\partial h_1}{\partial_{L_2}} & 0 & 0 & 0\\
0 & 0 & \frac{\partial h_3}{\partial_G} & \frac{\partial h_3}{\partial_{G_2}} & 0\\
0 & 0 & \frac{\partial^2 h_3}{\partial G^2}
      & \frac{\partial^2 h_3}{\partial G \partial {G_2}} & 0\\
0 & 0 & 0 & 0 & \frac{\partial h_4}{\partial J}
\end{pmatrix}.
\label{eq:HLY-rank-minor}
\end{equation}
The Keplerian factors are
\[
\frac{\partial h_0}{\partial L_1}
 = \frac{M_1^2\mu_1^3}{L_1^3} \neq0,
\qquad
\frac{\partial h_1}{\partial L_2} 
 = \frac{M_2^2\mu_2^3}{L_2^3} \neq0.
\]
For the $(G, G_2)$-block, it is therefore sufficient to consider
\begin{equation}
\Delta_3 =
\det
\begin{pmatrix}
\frac{\partial h_3}{\partial_G} & \frac{\partial h_3}{\partial_{G_2}}\\
\frac{\partial^2 h_3}{\partial G^2} & \frac{\partial^2 h_3}{\partial G \partial {G_2}}
\end{pmatrix}.
\label{eq:HLY-two-by-two-minor}
\end{equation}
In order to compute the determinant of $\Delta_3$ we must obtain the partial derivatives of $G_1^0$ with respect to $G$ and $G_2$ that we already needed when introducing the coordinates $x$,
$y$, see Appendix \ref{Appen}. We also need the second-order partials of $G_1^0$, namely, $\partial^2 G_1^0/\partial G^2$, $\partial^2 G_1^0/\partial G_2^2$ and $\partial^2 G_1^0/\partial G \partial {G_2}$, that are obtained
using \eqref{eq:G1equilibrium} and the first-order partials. Using the masses values of the 14 Herculis equilibrium and choosing $L_1$, $L_2$ and $G$ those of the equilibrium and $G_2$ as $G_2 = 3.15432 + \Delta G_2$, with $\Delta G_2$ in $[1/16, 1/12]$, and the corresponding $G_1^0$ annihilating the partial derivative of $\bar{H}_{\rm quad}$ with respect to $G_1$.
Substituting these values into $\Delta_3$ yields
$\det\Delta_3 \in [1.2847\times 10^{-7},\,1.3081\times 10^{-7}]$. Thus, $\det\Delta_3$ is positive throughout the considered interval of $G_2$.

\begin{figure*}
\centering
\captionsetup{width=\linewidth}
\includegraphics[width=1\textwidth]
{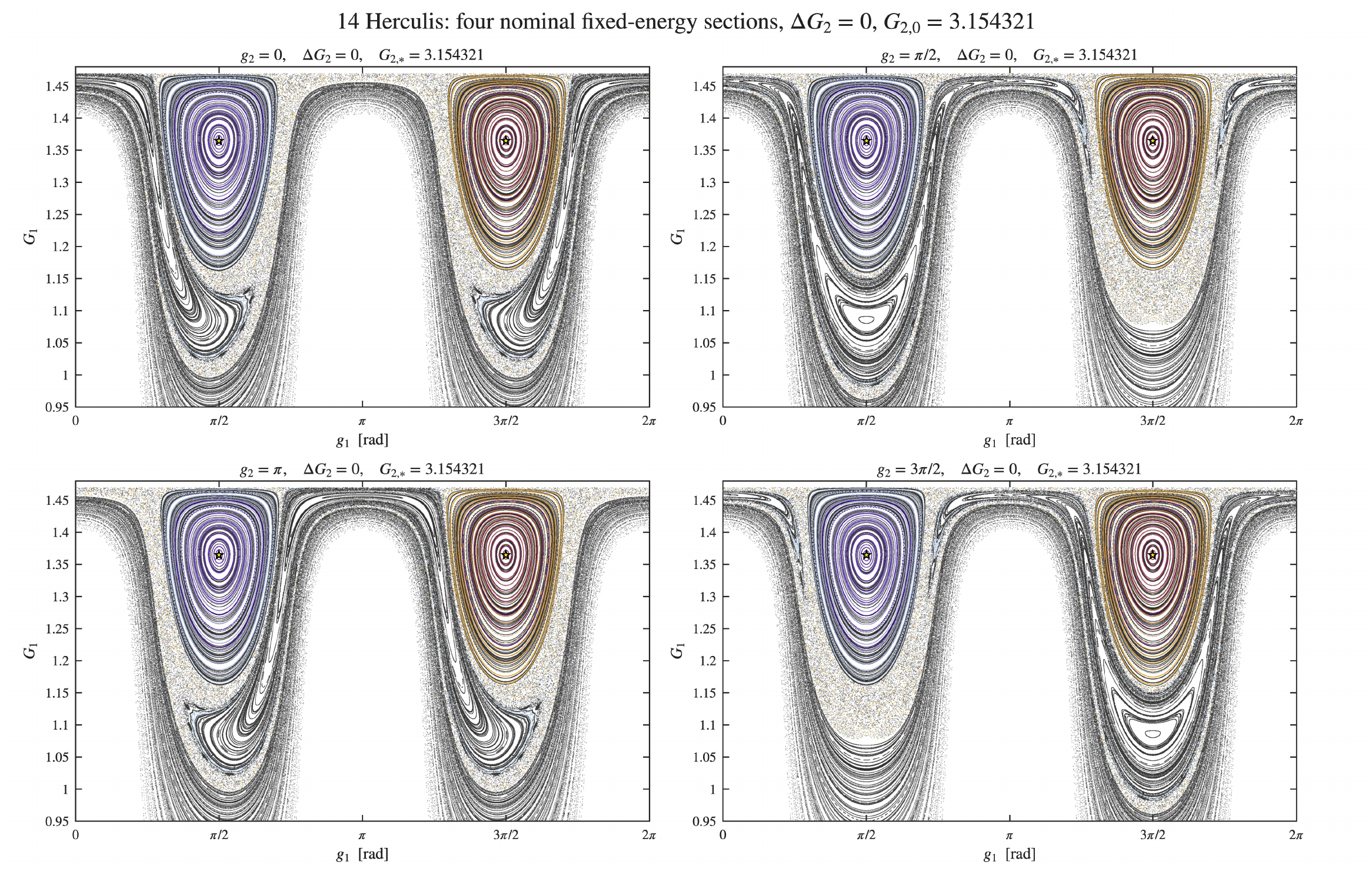}
\caption{
Fixed-energy Poincaré sections of the quadrupole--octupole secular
Hamiltonian for the nominal case \(\Delta G_2=0\), represented in
the \((g_1,G_1)\)-plane. The four panels correspond to
\(g_2=0,\pi/2,\pi,3\pi/2\), with \(\dot g_2>0\). The stars mark the
two elliptic relative equilibria of the quadrupole approximation.
The regular central islands persist throughout the sampled
\(g_2\)-cycle, while the outer regions exhibit a stronger
phase-dependent deformation.
}
\label{fig:poincare_qo_nominal}
\end{figure*}

Likewise, under the same assumptions on the actions as in the previous paragraph, the partial derivative \(\partial h_4/\partial J\) takes values approximately in the interval \([-0.00083, -0.00080]\). Moreover, we observe that \(\partial h_4/\partial J = 0\) if and only if \(\mathcal{Q}_b = 0\). However, this condition cannot be met when \(G > G_2\), subject to the constraints relating the actions \(L_1\), \(G\), \(G_2\), and \(G_s\). Consequently, we can conclude analytically, without resorting to numerical computations, that \(\partial h_4/\partial J\) does not vanish throughout the relevant domain we are considering.

Analogous to the previous paragraph, checking $\det \Delta_3$ numerically can be replaced by verifying that the partial derivative $\partial h_3 / \partial G_2$ is non-zero in our domain. Discarding the action $G$ then proves the existence of a full-rank $4 \times 4$ minor. This implies the persistence of invariant $4$-tori in an intermediate reduced space where $G$ is treated as an integral. Reintroducing $G$ as an action variable unreduces the symmetry, allowing us to reconstruct invariant $5$-tori.

As a result, we arrive at $\det\mathcal{R}_5 \neq 0$. By continuity, the same rank condition holds in an open neighbourhood of the reference configuration, and, by symmetry, for the branch associated with $g_1 = 3\pi/2$. 

Therefore, the hypotheses of Han-Li-Yi theorem \cite{HanLiYi2010} are satisfied when $g_1$, $G_1$ are in a neighbourhood of the equilibria associated with the 14 Herculis system. As a consequence, the invariant tori of the integrable approximation survive under the octupole and higher-order perturbations.

\medskip
\noindent
\textbf{Theorem 4.1.}
\emph{
For the Hamiltonian system associated with \eqref{eq:Hfull}, there exist two Cantor families of invariant KAM $5$-tori, associated with the two elliptic-type relative equilibria of the quadrupole secular approximation. The dynamics on these tori is quasi-periodic, and the corresponding trajectories remain confined to bounded regions of the twelve-dimensional phase space for all time. In particular, these quasi-periodic motions provide a dynamical scenario consistent with the observed configuration of the 14 Herculis system.
}

\medskip
\noindent
More precisely, the surviving KAM tori form a set of asymptotically full measure in the space of the action parameters $\mathcal{I}$. For a fixed local domain $D$ of positive measure, the Cantor set $D_\rho$ parametrising the persistent tori satisfies \cite{HanLiYi2010}
\[
\operatorname{meas}(D\setminus D_\rho)
=\mathcal{O}\!\left(\rho^{\sigma/s}\right),
\qquad 0<\sigma<\frac{1}{5},
\]
where $\rho = \varepsilon^{1/4}$ is the small parameter of the multiscale normal form and $s\geq 1$ is the derivative order entering the frequency-map non-degeneracy condition, thus $s=1$. (The choice of $\rho$ is motivated by the need to properly account for the natural powers of the small parameter; in view of \eqref{HamNor}, this leads to the appropriate scaling.) As a consequence the relative measure of the excluded parameters tends to zero as $\rho \to 0$. This estimate bounds the set for which persistence is not guaranteed, without implying that all corresponding tori are destroyed.

This result provides a rigorous explanation for the long-term stability of the observed architecture. Rather than being organised by isolated periodic orbits, the dynamics is structured by two families of invariant quasi-periodic motions that surround the relative equilibrium. The observed system may therefore be interpreted as evolving on, or sufficiently close to, one of these invariant tori.

When $G_2$ is restricted to the interval $[3.2168, 3.2376]$, the corresponding value of $G_1^0$ ranges 
from $1.3092$ to $1.3233$. Consequently, for the quadrupolar relative equilibrium, the orbital eccentricities 
vary approximately between $0.4348$ and $0.4542$ for the inner planet, and between $0.6148$ and 
$0.6212$ for the outer planet. For the associated quasi-periodic motions, however, $G_1$ and $G_2$ are allowed 
to vary—-specifically by about $\pm 0.05$ for $G_1$ and $\pm 0.02$ for $G_2$ around their 
equilibrium values. Under these conditions, the eccentricities oscillate between $0.356$ and $0.515$ 
for 14 Her b, and between $0.610$ and $0.627$ for 14 Her c, remaining fully compatible with reported orbital 
parameters. Furthermore, using these estimations for $G_1$ and $G_2$, the mutual inclination between the two 
planetary orbital planes ranges from $45.440^\circ$ to $56.366^\circ$, which is also consistent with 
observation data.

%%%%%%%%%%%%%%%%%%%%%%%%%%%%%%%%%%%%%%%%%%%%%%%%

\section{Numerical experiments}
\label{sec:numerical}

\begin{figure*}
\centering
\captionsetup{width=\linewidth}
\includegraphics[width=1\textwidth]
{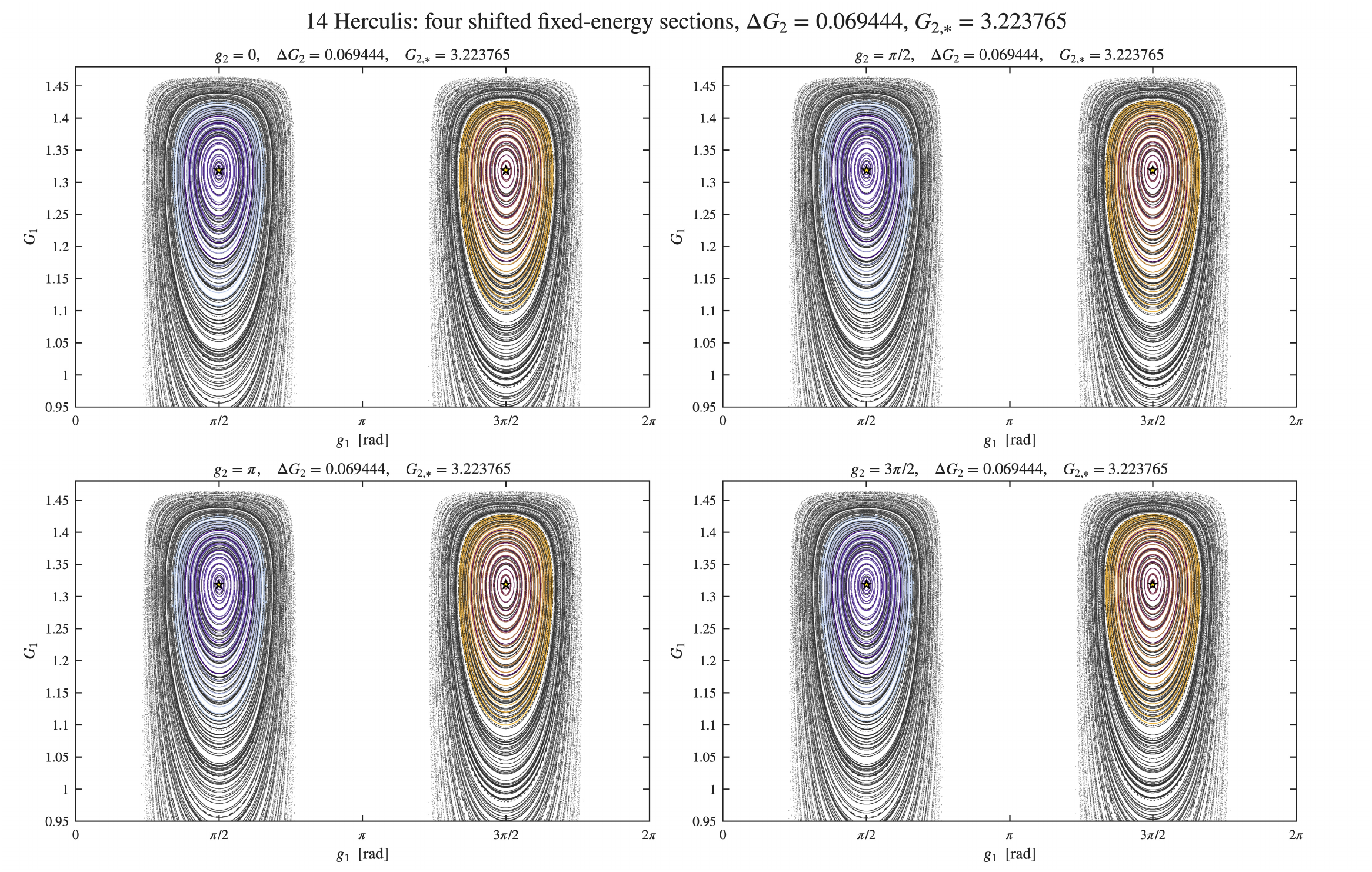}
\caption{
Same as Figure~\ref{fig:poincare_qo_nominal}, but for the displaced
case \(\Delta G_2=5/72\), with
\(G_{2,\star}=3.223765\). The two local elliptic islands remain
clearly organised throughout the \(g_2\)-cycle. In contrast with the
nominal case, the outer accessible regions form two visually
separated components.
}
\label{fig:poincare_qo_shifted}
\end{figure*}

In contrast with the integrable quadrupole model, the inclusion of
the octupole contribution produces a secular Hamiltonian with two
degrees of freedom. We therefore examine its dynamics by means of
two-dimensional Poincaré sections on fixed energy levels.

In addition to the nominal configuration, we numerically examine the
interval of outer angular-momentum actions considered in
Section~4.2. More precisely, we write
\begin{equation}
G_{2,\star}=G_{2,0}+\Delta G_2,
\qquad
G_{2,0}=3.15432,
\label{eq:numerical_G2_shift}
\end{equation}
and consider
\begin{equation}
\Delta G_2\in
\left\{
0,\frac{1}{16},\frac{5}{72},
\frac{11}{144},\frac{1}{12}
\right\},
\label{eq:numerical_G2_shifts}
\end{equation}
where the four non-zero values are equally spaced
samples of the interval $\Delta G_2\in\left[{1}/{16},{1}/{12}\right]$ used to arrange the Hamiltonian in the multiscale regime required for the KAM analysis. For every value of \(G_{2,\star}\), the corresponding quadrupolar equilibrium
\[
G_{1,\star}=G_1^0(G_{2,\star})
\]
is recomputed from the equilibrium equation \eqref{eq:G1equilibrium}. The energy level used
in each numerical experiment is then defined consistently by
\begin{equation}
E_{\Delta G_2}^{\rm qo}
=
\overline H_{\rm qo}
\left(
\frac{\pi}{2},G_{1,\star},0,G_{2,\star}
\right)
=
\overline H_{\rm qo}
\left(
\frac{3\pi}{2},G_{1,\star},0,G_{2,\star}
\right),
\label{eq:shifted_energy_qo}
\end{equation}
where
\[
\overline H_{\rm qo}
=
\overline H_{\rm quad}
+
\overline H_{\rm octu}.
\]
Thus, changing \(G_{2,\star}\) also changes \(G_{1,\star}\) and the
energy surface. The nominal and displaced sections should therefore
not be interpreted as different samplings of a common energy level.

We first consider oriented sections of the form
\begin{equation}
g_2=\gamma,
\qquad
\dot g_2>0,
\qquad
\overline H_{\rm qo}=E_{\Delta G_2}^{\rm qo},
\label{eq:g2_sections_shifted}
\end{equation}
and represent their intersections in the \((g_1,G_1)\)-plane. To
sample the principal phases of $g_2$, we produce sections at \(g_2=0\) and \(g_2=\pi\), where the contributions proportional to \(\sin g_2\) vanish. The
sections at \(g_2=\pi/2\) and \(g_2=3\pi/2\), for which
\(\lvert\sin g_2\rvert=1\), test the dynamics where this additional
local suppression is absent.

Figure~\ref{fig:poincare_qo_nominal} shows the nominal case
\(\Delta G_2=0\). The two elliptic islands surrounding
\(g_1=\pi/2\) and \(g_1=3\pi/2\) remain clearly identifiable in all
four sections. Their interiors are foliated by nested invariant
curves, whereas the outer parts of the section display a stronger
dependence on \(g_2\) and contain irregularly populated layers. On
the nominal energy surface, these outer layers extend through the
large-amplitude part of the section and visually connect the
neighbourhoods of the two elliptic regions.

\begin{figure*}
\centering
\captionsetup{width=\linewidth}
\includegraphics[width=0.96\textwidth]
{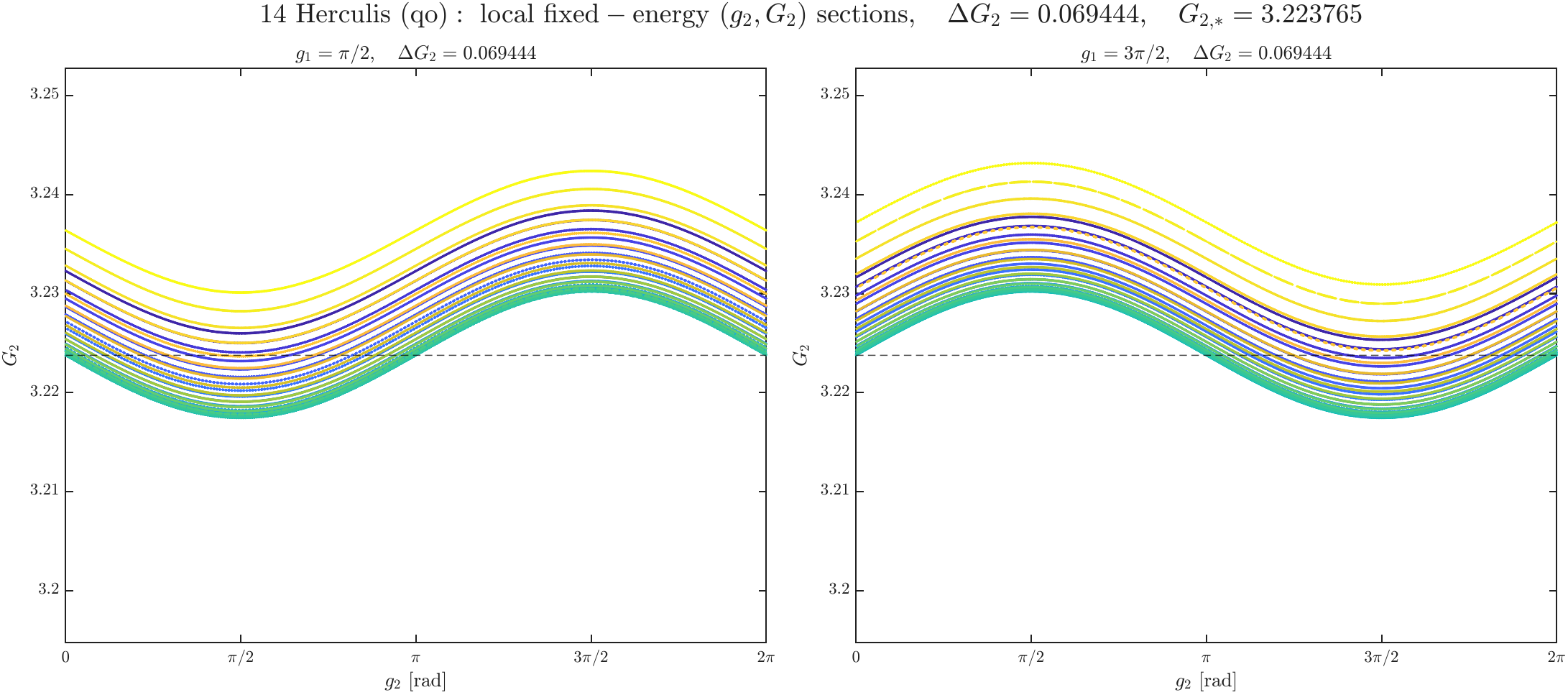}
\caption{
Fixed-energy Poincaré sections of the quadrupole--octupole secular
Hamiltonian for \(\Delta G_2=5/72\), represented in the
\((g_2,G_2)\)-plane. The left and right panels correspond to
\(g_1=\pi/2\) and \(g_1=3\pi/2\), respectively, with
\(\dot g_1>0\). The curves wind around the complete angular
direction, demonstrating circulation of \(g_2\), while \(G_2\)
undergoes a small, smooth modulation.
}
\label{fig:poincare_g2G2_shifted}
\end{figure*}

The behaviour changes noticeably after displacing \(G_2\). As a
representative example, Figure~\ref{fig:poincare_qo_shifted} shows
the case
\[
\Delta G_2=\frac{5}{72},
\qquad
G_{2,\star}=3.223765,
\qquad
G_{1,\star}=1.318645.
\]
This value is the second of the four equally spaced displacements in
\([1/16,1/12]\). The other three
non-zero displacements give the same qualitative picture.

The local elliptic organisation survives the displacement: both
islands remain surrounded by smooth nested curves in all four
sections. The global geometry, however, is reorganised. The broad
outer layer present in the nominal case is strongly reduced, and the
accessible regions surrounding the two islands appear as two
separate components. This difference is a consequence of following
the displaced quadrupolar equilibrium to the new energy
\(E_{\Delta G_2}^{\rm qo}\), rather than merely changing the initial
conditions on the nominal surface.

To verify that the missing connection in
Figure~\ref{fig:poincare_qo_shifted} was not caused by an
insufficient number or an inappropriate placement of initial
conditions, we performed an auxiliary fixed-energy accessibility
scan. For prescribed values of \((g_1,G_1)\) and of the section angle
\(\gamma\), we searched all physically admissible \(G_2\) satisfying
\begin{equation}
F_{\gamma,g_1,G_1}(G_2)
=
\overline H_{\rm qo}(g_1,G_1,\gamma,G_2)
-
E_{5/72}^{\rm qo}
=0.
\label{eq:energy_accessibility_equation}
\end{equation}
We used the 32 section angles
\[
\gamma_j=\frac{j\pi}{16},
\qquad
j=0,\ldots,31,
\]
and examined both possible passages between the islands: the central
gate near \(g_1=\pi\) and the periodic gate near
\(g_1=0\equiv2\pi\). The action interval was extended to
\[
0.10\leq G_1\leq L_1-0.002,
\]
and the complete admissible interval of \(G_2\) was searched,
including possible tangent roots that do not produce a sign change.

No root of \eqref{eq:energy_accessibility_equation} was detected in
either gate for any of the sampled section angles. Moreover, the
smallest energy mismatch was
\[
\min
\left|
\overline H_{\rm qo}
-
E_{5/72}^{\rm qo}
\right|
\simeq 2.44\times10^{-2},
\]
which is well separated from the numerical root tolerance. Within
the explored domain, this indicates that the absence of a bridge
between the two islands is a property of the displaced energy
surface rather than a sampling artefact. The very small interval
\(0<G_1<0.10\), not included in this test, corresponds to an almost
rectilinear inner orbit and lies outside the neighbourhood relevant
to the elliptic configurations studied here.

A complementary view is obtained by reversing the roles of the
secular angles. We consider the oriented sections
\begin{equation}
g_1=\frac{\pi}{2},
\qquad
\dot g_1>0,
\qquad
\overline H_{\rm qo}=E_{\Delta G_2}^{\rm qo},
\label{eq:g1_pi2_section}
\end{equation}
and
\begin{equation}
g_1=\frac{3\pi}{2},
\qquad
\dot g_1>0,
\qquad
\overline H_{\rm qo}=E_{\Delta G_2}^{\rm qo},
\label{eq:g1_3pi2_section}
\end{equation}
and represent their intersections in the \((g_2,G_2)\)-plane.  Figure~\ref{fig:poincare_g2G2_shifted} shows the representative
displaced case \(\Delta G_2=5/72\). The curves extend throughout the
complete interval \(0\leq g_2<2\pi\), showing that \(g_2\) circulates
rather than librates around \(0\) or \(\pi\). Its conjugate action
\(G_2\) undergoes only a small and smooth modulation. The two
branches have approximately opposite modulation phases, consistently
with the symmetry between \(g_1=\pi/2\) and \(g_1=3\pi/2\). The same
qualitative behaviour is obtained for the nominal case and for the
other three displaced values.

Finally, we assess the robustness of these conclusions with respect
to the next multipole of the secular expansion. We introduce the
quadrupole--octupole--hexadecapole Hamiltonian
\begin{equation}
\overline H_{\rm qoh}
=
\overline H_{\rm quad}
+
\overline H_{\rm octu}
+
\overline H_{\rm hexa},
\qquad
\overline H_{\rm hexa}
=
\mathcal{O}(\delta^5),
\label{eq:qoh_hamiltonian}
\end{equation}
where $\overline H_{\rm hexa}$ is given in the Appendix~\ref{Appen:Hexadecapole}. Next, we repeat the complete numerical experiment. For each value of
\(\Delta G_2\), the energy is recomputed as
\begin{equation}
E_{\Delta G_2}^{\rm qoh}
=
\overline H_{\rm qoh}
\left(
\frac{\pi}{2},G_{1,\star},0,G_{2,\star}
\right),
\label{eq:shifted_energy_qoh}
\end{equation}
and the initial conditions are reconstructed on the corresponding
fixed-energy surface.

\begin{figure*}
\captionsetup{width=\linewidth}
\centering
\includegraphics[width=0.495\textwidth]
{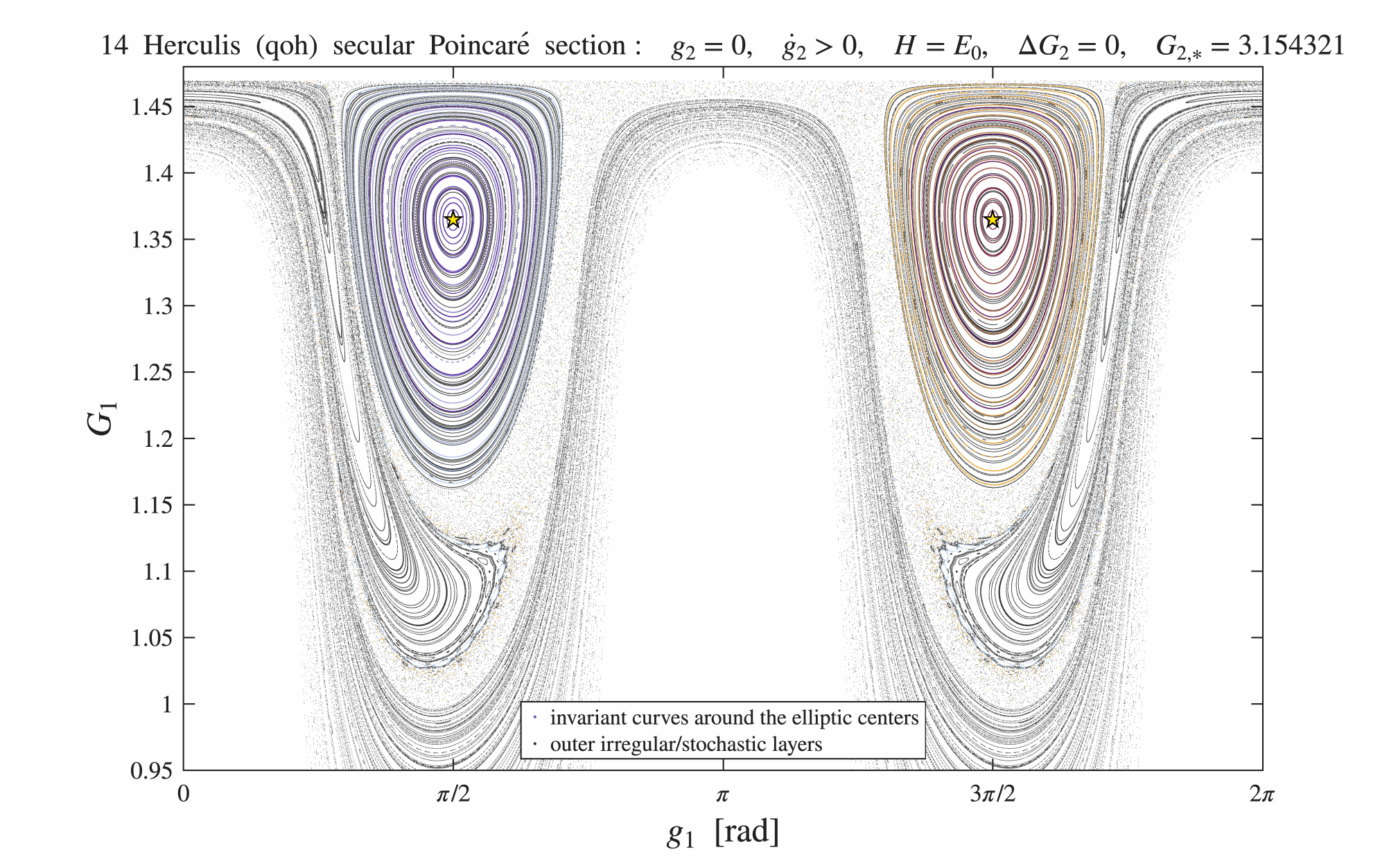}
\includegraphics[width=0.495\textwidth]
{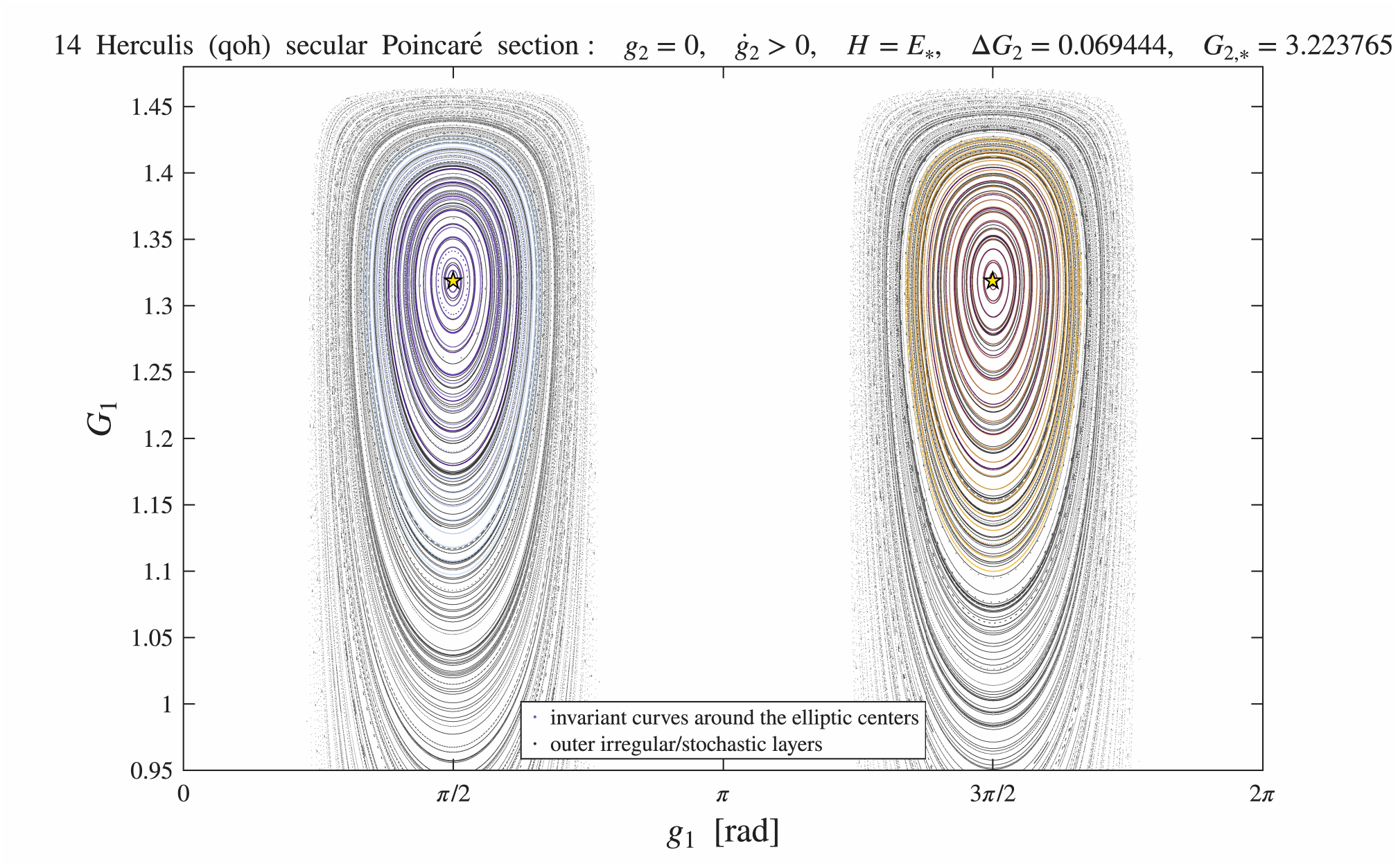}
\caption{
Fixed-energy Poincaré sections of the
quadrupole--octupole--hexadecapole secular Hamiltonian for
\(\Delta G_2=0\) and \(\Delta G_2=5/72\), represented in the \((g_1,G_1)\)-plane at
\(g_2=0\), with \(\dot g_2>0\). Comparison with Figure~\ref{fig:poincare_qo_nominal} and 
Figure~\ref{fig:poincare_qo_shifted} shows no qualitative change in
the regular elliptic islands or in the surrounding global
structure.
}
\label{fig:poincare_qoh_shifted}
\end{figure*}

At the five reference configurations, the relative hexadecapole
contribution satisfies approximately
\[
8\times10^{-4}
\lesssim
\frac{
\left|\overline H_{\rm hexa}\right|
}{
\left|\overline H_{\rm quad}+\overline H_{\rm octu}\right|
}
\lesssim
2.3\times10^{-3}.
\]
The corresponding sections are almost indistinguishable from those
of the quadrupole--octupole truncation. As an illustration,
Figure~\ref{fig:poincare_qoh_shifted} shows the same displaced case
as Figure~\ref{fig:poincare_qo_nominal} and Figure~\ref{fig:poincare_qo_shifted} for $g_2=0$, now including the
hexadecapole term. The locations and sizes of the elliptic islands,
the deformation of the outer layers, and the separation between the
two accessible components remain unchanged at the qualitative
level.

Taken together, these experiments show that the elliptic structures
organised by the quadrupole approximation persist throughout the
circulation of \(g_2\), over the complete displaced interval
\(\Delta G_2\in[1/16,1/12]\), and after inclusion of the
hexadecapole correction. The displacement of \(G_2\) does modify the
global geometry of the fixed-energy section, most visibly by
separating the two accessible components, but it does not destroy the
local families of regular motions surrounding either elliptic
branch.

%%%%%%%%%%%%%%%%%%%%%%%%%%%%%%%%%%%%%%%%%%%%%%%

\section{Conclusions}
\label{sec:conc}

We have developed a Hamiltonian interpretation of the long-term
dynamics of the 14 Herculis planetary system. Starting from the
spatial planetary three-body problem, we have employed canonical
heliocentric variables, performed the symplectic reduction of the
nodes using coordinates adapted to the total angular momentum, and  
averaged over the two fast Keplerian angles. The hierarchy between
the planetary masses and the stellar mass, together with the
separation between the two orbital scales, allows us to
organise the secular Hamiltonian into quadrupolar, octupolar, and
higher-order contributions. This combination of geometric
reduction, suitable canonical coordinates, averaging theory, and KAM theory
within the Hamiltonian perturbation setting, provides a mathematically
controlled framework in which the eccentric and mutually inclined
architecture of 14 Herculis can be studied.

For parameters compatible with the observationally inferred
configuration, the quadrupole secular Hamiltonian possesses an
elliptic relative equilibrium associated with the branches
\(g_1=\pi/2\) and \(g_1=3\pi/2\). The equilibria obtained here are
not intended to reproduce the nominal orbital fit exactly, rather,
they act as organising centres for 
families of nearby secular motions. When the octupole contribution is
incorporated, the integrability of the quadrupole approximation is
lost and a second secular degree of freedom arises. Nevertheless,
near \(g_1 = \pi/2\) or \(3\pi/2\) and $G_1$ near $G_1^0$, the secular Hamiltonian retains the
properly degenerate multiscale structure required by 
Han-Li-Yi theorem \cite{HanLiYi2010}. Under the corresponding
smallness, non-resonance, and non-degeneracy conditions, the
periodic orbit issued from the quadrupolar relative equilibrium is
surrounded by a family of invariant \(2\)-tori in the reduced
secular space \(\mathcal{S}_{L_1,L_2,G}\).

Reintroducing the two fast mean anomalies eliminated by the
averaging procedure, together with the cyclic angle associated with
the angular-momentum reduction, reconstructs these invariant
structures as five-dimensional KAM tori of the original spatial
planetary problem. The quasi-periodic motions carried by these tori
remain confined for all time. Our result therefore supplies a
rigorous Hamiltonian mechanism for the long-term confinement of
configurations compatible with 14 Herculis. 

By slightly modifying the range of variation of \(G_2\) from that inferred from the reported eccentricities, 
we can separate the different orders of the principal part of the Hamiltonian, as required for the application 
of KAM theory. Under these conditions, the eccentricity of 14 Her b varies approximately between \(0.356\) 
and \(0.515\), while that of 14 Her c ranges from \(0.610\) to \(0.627\). Both intervals remain consistent 
with the reported orbital parameters. Furthermore, the corresponding estimates of the mutual inclination 
between the two planetary orbital planes range from \(45.440^\circ\) to \(56.366^\circ\), also in agreement 
with observational data.

From the perpestive of KAM theory, the architecture of 14 Herculis, while
observationally remarkable, is not dynamically anomalous within the
spatial three-body problem. Large eccentricities and substantial
mutual inclinations do not necessarily imply instability; they may
instead belong to stable regions organised by elliptic relative
equilibria and their surrounding KAM tori. The method presented here is not tied to
one particular orbital fit or to 14 Herculis itself. The same
strategy may be applied to other hierarchical planetary systems.

For one star and two planets, the relation between elliptic relative
equilibria and five-dimensional invariant tori in the spatial
three-body problem is supported by a substantial theoretical
foundation \cite{PalacianSayasYanguas2013,
PalacianSayasYanguas2015}. The corresponding problem for a star
with three or more planets is considerably less developed,
especially for eccentric and mutually inclined configurations that
are not perturbations of circular and coplanar motions. In such
systems one expects higher-dimensional invariant tori of maximal dimension, at least tori of 
dimension eight when the system is formed by a star and three planets. Linearly stable relative equilibria in suitable
reduced secular spaces provide promising candidates from which such
KAM families could be constructed. This opens a route towards
analytical studies of multi-planet extrasolar systems and nested
planet-satellite hierarchies.

At the same time, stability should not be interpreted as a global
property of the planetary phase space. 
%Complementary work in the
%planetary four-body problem has %established the existence of
%trajectories displaying large %variations of the orbital %elements
%and instability of the semi-%major axes
%\cite{ClarkeFejozGuardia2024}. 
The coexistence of instability
channels with Cantor families of invariant KAM tori highlights the
richness of higher-dimensional planetary dynamics. A central
challenge for future work is therefore not merely to classify a
system as stable or unstable, but to identify the 
organising structures that determine which regions of its phase
space support confined quasi-periodic motion and which permit
large-scale transport.

An additional direction for future work would be to retain the rotational dynamics and finite-size effects of planets or satellites, building on the Poisson and invariant-based reduction frameworks developed for the full gravitational \(N\)-body problem \cite{CrespoTurner2022,CrespoDullin2026}.

%%%%%%%%%%%%%%%%% APPENDICES %%%%%%%%%%%%%%%%%%%%%

\appendix

\section{Completing the transformations}
\label{Appen}

The variables introduced in Section \ref{QuasiKAM} require slight modifications of some of their angle coordinates in order for the transformations to be canonical.

We begin with the transformation from the set $\{ l_1, $ $l_2, $ $g_1, $ $g_2, $ $\nu_1, $ $\nu, $ $L_1, $ $L_2,$ $G_1, $ $G_2,$ $G,$ $ B\}$ to
$\{ \bar{l}_1, $ $l_2, $ $x, $ $\bar{g}_2, $ $\bar{\nu}_1, $ $\nu, $ $L_1, $ $L_2, $ $y, $ $G_2, $ $G, $ $B \}$ given
by 
\[ x = g_1 - \pi/2\quad \mbox{or} \quad x = g_1 - 3\pi/2, \quad y = G_1 - G_1^0. \] 
As $G_1^0$ depends on $L_1$, $G_2$ and $G$, the associated transformed angles become
\[ 
\bar{l}_1 = l_1 + g_1 \frac{\partial G_1^0}{\partial L_1}, \quad 
\bar{g}_2 = g_2 + g_1 \frac{\partial G_1^0}{\partial G_2}, \quad 
\bar{\nu}_1 = {\nu}_1 + g_1 \frac{\partial G_1^0}{\partial G}, \]
where the partial derivatives of $G_1^0$ with respect to the actions are obtained from \eqref{eq:G1equilibrium}. It yields
\[
\begin{aligned}
\frac{\partial G_1^0}{\partial L_1} =\,& \frac{5L_1 \left((G^2 - G_2^2)^2 - G_s^2 \right)}{2 G_s^{3/2} (5L_1^2 + 4G_2^2 + 8 G^2 - 12 G_s)}, \\
\frac{\partial G_1^0}{\partial G_2} =\,& -\frac{G_2 \left( 5L_1^2 (G^2 - G_2^2) + 2G_s^2 \right)}{G_s^{3/2} (5L_1^2 + 4G_2^2 + 8 G^2 - 12 G_s)}, \\
\frac{\partial G_1^0}{\partial G} =\,& \frac{5L_1^2 G (G^2-G_2^2) - 4 G G_s^2}{G_s^{3/2} (5L_1^2 + 4G_2^2 + 8 G^2 - 12 G_s)}.
\end{aligned}    
\]
The denominators take positive values in the region of phase space where we perform our analysis.

%%%%%%%%%%%%%%%%%%%%%%%%%%%%%%%%%%%%%%%%%%%%%%%%%%
The second transformation introduces the set $\{ \tilde{l}_1, $ $l_2, $ $\varphi, $ $\tilde{g}_2, $ $\tilde{\nu}_1, $ $\nu, $ $L_1, $ $L_2, $ $J, $ $G_2, $ $G, $ $B \}$
from $\{ \bar{l}_1, $ $l_2, $ $x, $ $\bar{g}_2, $ $\bar{\nu}_1, $ $\nu, $ $L_1, $ $L_2, $ $y, $ $G_2, $ $G, $ $B \}$ by means of
\[ \varphi = \mbox{arctan2}(y, \omega x), \quad J = (\omega^2 x^2 + y^2)/(2\omega).\]
The frequency $\omega$ has been given in \eqref{freq}, and as it depends on $L_1$, $G_2$ and $G$, their conjugate angles $\bar{l}_1$, $\bar{g}_2$ and $\bar{\nu}_1$ are corrected
through
\[ 
\tilde{l}_1 = \bar{l}_1 + \frac{x y}{2\omega}\frac{\partial \omega}{\partial L_1}, \quad 
\tilde{g}_2 = \bar{g}_2 + \frac{x y}{2\omega}\frac{\partial \omega}{\partial G_2}, \quad 
\tilde{\nu}_1 = \bar{\nu}_1 + \frac{x y}{2\omega}\frac{\partial \omega}{\partial G}.\]
The partial derivatives of $\omega$ with respect to $L_1$, $G_2$, and $G$ are obtained directly 
from \eqref{freq}; their explicit expressions are omitted here for brevity, as they are rather lengthy.

\section{The double-averaged hexadecapole contribution}
\label{Appen:Hexadecapole}

This section contains the hexadecapole term contribution, which can be written as
\begin{equation}
\overline H_{\rm hexa}
=
\kappa'' h_{\rm hexa},
\label{eq:Hhexa}
\end{equation}
where
\begin{equation}
\kappa''
=
-\frac{
9m_1m_2M_2^5\mu_2^{10}L_1^8
}{
8192M_1^4\mu_1^8L_2^{10}
}.
\label{eq:kappa-hexa}
\end{equation}
To express the result compactly, we use
$\chi = \cos(I_1 + I_2)$, which was given in (\ref{eq:mutualinclination}).
The dimensionless hexadecapole Hamiltonian is then
\begin{align}
h_{\rm hexa}
={}&
\frac{1}{\eta_2^7}
\Big[
 \mathcal A_0
+\mathcal A_1\cos(2g_1)
+\mathcal A_2\cos(2g_2)
\nonumber\\
&\quad
+\mathcal A_{+}\cos\bigl(2g_1+2g_2\bigr)
+\mathcal A_{-}\cos\bigl(2g_1-2g_2\bigr)
\nonumber\\
&\quad
+\mathcal B_0\cos(4g_1)
+\mathcal B_{+}\cos\bigl(4g_1+2g_2\bigr)
\nonumber\\
&\quad
+\mathcal B_{-}\cos\bigl(4g_1-2g_2\bigr)
\Big],
\label{eq:hhexa}
\end{align}
where
\begin{align}
\mathcal A_0
={}&
(2+3e_2^2)
(8+40e_1^2+15e_1^4)
(35\chi^4-30\chi^2+3),
\label{eq:hexa-A0}
\\[1mm]
\mathcal A_1
={}&
140e_1^2(2+e_1^2)(2+3e_2^2)
(1-\chi^2)(7\chi^2-1),
\label{eq:hexa-A1}
\\[1mm]
\mathcal A_2
={}&
10e_2^2(8+40e_1^2+15e_1^4)
(1-\chi^2)(7\chi^2-1),
\label{eq:hexa-A2}
\\[1mm]
\mathcal A_{+}
={}&
140e_1^2e_2^2(2+e_1^2)
(1-\chi)^2
(7\chi^2+7\chi+1),
\label{eq:hexa-Aplus}
\\[1mm]
\mathcal A_{-}
={}&
140e_1^2e_2^2(2+e_1^2)
(1+\chi)^2
(7\chi^2-7\chi+1),
\label{eq:hexa-Aminus}
\\[1mm]
\mathcal B_0
={}&
735e_1^4(2+3e_2^2)
(1-\chi^2)^2,
\label{eq:hexa-B0}
\\[1mm]
\mathcal B_{+}
={}&
735e_1^4e_2^2
(1-\chi)^3(1+\chi),
\label{eq:hexa-Bplus}
\\[1mm]
\mathcal B_{-}
={}&
735e_1^4e_2^2
(1-\chi)(1+\chi)^3.
\label{eq:hexa-Bminus}
\end{align}

Consequently, the secular Hamiltonian retained for the higher-order
numerical comparison is
\begin{equation}
H_{\rm sec}^{\rm qoh}
=
H_{\rm Kep}
+
\varepsilon
\left(
\overline H_{\rm mono}
+
\overline H_{\rm quad}
+
\overline H_{\rm octu}
+
\overline H_{\rm hexa}
\right).
\label{eq:Hsec-qoh}
\end{equation}

\section*{Acknowledgements}
The author ET gratefully acknowledges the support of Universidad Viña del Mar through the research project FIIUVM-CYT-2609. The authors JFP and PY acknowledge the support of Project PID2022-140469NB-C21 of Fondo Europeo de Desarrollo Regional, Agencia Estatal de Investigaci\'on (AEI) of Spain.

%%%%%%%%%%%%%%%%%%%%%%%%%%%%%%%%%%%%%%%%%%%%%%%%%%
\section*{Data Availability}
No data were used in the preparation of this manuscript
%%%%%%%%%%%%%%%%%%%% REFERENCES %%%%%%%%%%%%%%%%%%

\bibliographystyle{abbrv}
\bibliography{example}

\end{document}